\documentclass[12pt]{amsart}
\usepackage{pst-tree}

\newtheorem{thm}{Theorem}[section]
\newtheorem{thrm}[thm]{Theorem}

\newtheorem{prop}[thm]{Proposition}
\newtheorem{coro}[thm]{Corollary}

\theoremstyle{definition}
\newtheorem{defn}[thm]{Definition}
\newtheorem{remk}[thm]{Remark}
\newtheorem*{ackn}{Acknoledgements}

\newcommand{\1}{\operatorname{\mathbf{1}}}
\newcommand{\coDG}{\operatorname{{\it col}\mathcal{DG}_\infty^\bullet}}
\newcommand{\mDG}{\operatorname{{\it m}\mathcal{DG}_\infty^\bullet}}

\begin{document}

\title[Algebraic String Operations]
{Algebraic String Operations}

\author[T.~Tradler]{Thomas~Tradler}
\address{Thomas Tradler, Department of Mathematics, College of Technology of the City University
of New York, 300 Jay Street, Brooklyn, NY 11201, USA}
\email{ttradler@citytech.cuny.edu}

\author[M.~Zeinalian]{Mahmoud~Zeinalian}
\address{Mahmoud Zeinalian, Department of Mathematics, C.W. Post Campus
of Long Island University, 720 Northern Boulevard, Brookville, NY
11548, USA} \email{mzeinalian@liu.edu}

\begin{abstract}
We first discuss how open/closed chord diagrams, both with and
without marked points, act on appropriate Hochschild complexes
possibly coupled with the two-sided cobar complex. Then, in the
main part of the paper, we introduce the notion of a V$_k$-algebra
and obtain suitable homotopy versions.
\end{abstract}

\maketitle

\setcounter{tocdepth}{1}
\tableofcontents

\section{Introduction}\label{intro}

It is long known that the Hochschild cohomology of an associative
algebra has the structure of a Gerstenhaber algebra; see \cite{G}.
In recent years, it has become clear that this is only a starting
point for the study of the algebraic structure of the Hochschild
complex. Deligne's conjecture states that the operad of chains on
the little disc operad acts on the Hochschild cochain complex of
an associative algebra. In the more general case of
A$_\infty$-algebras, this conjecture was proved in \cite{KoSo}.
The original associative case was also proved in \cite{Ka1},
\cite{MS1}, \cite{Ta1}, \cite{Ta2}, \cite{V}. The discovery of
String Topology \cite{CS1} suggested that in the presence of an
appropriate notion of Poincar\'e duality a richer collection of
operations should exist. In fact, under suitable assumptions,
there is a plethora of operations both in the Hochschild and
cyclic Hochschild settings.

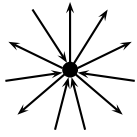
\begin{figure}
\[
\begin{pspicture}(0,0)(2,2)
 \psdots[dotsize=6pt](1,1)
 \psline[arrows=->](1,1)(0.3,0.4)         \psline[arrows=->](1,1)(1,1.9)
 \psline[arrows=->](1,1)(1.7,0.4)
 \psline[arrows=<-](1.82,1.37)(1.1,1.02)  \psline[arrows=<-](0.18,1.37)(0.9,1.02)
 \psline[arrows=<-](1.5,1.8)(1,1)         \psline[arrows=->](0.5,1.8)(0.96,1.09)
 \psline[arrows=->](0.8,0.2)(0.98,0.9)    \psline[arrows=->](1.2,0.2)(1.02,0.9)
 \psline[arrows=->](0.14,0.85)(0.9,0.96)  \psline[arrows=->](1.86,0.85)(1.1,0.96)
\end{pspicture}
\] \caption{The general vertex}
 \label{general-vertex}
\end{figure}

A step in this direction was taken in our previous paper \cite{TZ},
where the Hochschild complex of an associative algebra with an
invariant and symmetric co-inner product was treated. We will recall
the action of the cyclic Sullivan chord diagrams on the Hochschild
complex of an associative algebra with values in its dual from
\cite{TZ} and extend the discussion to the action of the chord
diagrams without marked points on the cyclic Hochschild complex, as
well as the action of open/closed chord diagrams on the cyclic
Hochschild complex coupled with the two-sided cobar complex. The
main part of this paper is a study of the algebraic structure of the
cyclic Hochschild complex of an A$_\infty$-algebra with an
appropriate homotopy co-inner product. An appropriate homotopy
version of an invariant and symmetric co-inner product is formulated
using the notion of a V$_\infty$-algebra. This consists of a system
of elements labelled by vertices with cyclically ordered directed
edges, as shown in Figure \ref{general-vertex}.
For each such vertex, with $m$ incoming edges shuffled in between
$n$ outgoing edges according to a partition $m=i_1 + \cdots +i_n$,
we associate an element of $A \otimes (A^*)^{\otimes i_1} \otimes
\cdots \otimes A \otimes (A^*)^{\otimes i_n}$, where $A$ is the
underlying vector space of our algebra, and $A^*$ its dual. These
elements must satisfy certain conditions, which are stated in
Definition \ref{V_k}. If we restrict the structure to vertices with
only one outgoing edge, then we recover the notion of an
A$_\infty$-algebra. In general, for a fixed $k=1,2,\cdots, \infty$,
we consider V$_k$-algebras, governed by of vertices with $m$
incoming and $n$ outgoing edges, where $m$ is arbitrary while $n <
k+1$. In general, there are obstructions to completing a
V$_k$-algebra to a $V_{k+1}$-algebra. Interesting examples of
V$_k$-algebras occur in topology. For instance, $V_2$-algebras,
which are A$_\infty$-algebras with an invariant and symmetric
homotopy co-inner product, can be constructed from a triangulated
and oriented Poincar\'e duality space; see \cite{Z}.

We will describe  graphs which generalize cyclic Sullivan Chord
diagrams and construct a graph complex, denoted by
$\mathcal{DG}_k^\bullet$.  Roughly speaking,
$\mathcal{DG}_k^\bullet$ is generated by the set of directed
graphs with vertices having at most $k$ outgoing edges and has the
structure of a PROP. Equipped with these notions, we show how the
cyclic Hochschild complex of a V$_k$-algebra has a natural action
of $\mathcal{DG}_k^\bullet$, for $k=1,2,\cdots, \infty$. The case
of $k=1$ is somewhat akin to what was considered in \cite{KoSo}
for the Hochschild complex of an A$_\infty$-algebra. In this
paper, we treat the cyclic Hochschild complex in detail because it
is very clean and the graphs do not carry any extra decorations
such as marked points. We will then discuss modifications needed
to address the open/closed string interactions, i.e., operation on
the cyclic Hochschild complex coupled with the two-sided cobar
complex, as well as the case in which the cyclic Hochschild
complex is replaced by the Hochschild complex.

Here is a short description of the individual sections. Section
\ref{strict-chapter} is devoted to strictly associative algebras
with an invariant and symmetric co-inner products. The open/closed
case and the cases of chord diagrams both with and without marked
points will be investigated. In section \ref{algebra section}, the
notion of V$_k$-algebras will be defined and several special cases
will be examined. In section \ref{action section}, we define the
PROP $\mathcal{DG}_k^\bullet$ and show that it acts on the cyclic
Hochschild complex of a V$_k$-algebra. Sections \ref{algebra
section} and \ref{action section} form the core of this paper. In
section \ref{variations section}, we discuss the modifications
necessary to obtain the open/closed and the non-equivariant
operations for homotopy algebras. Also, the relevance of the
algebraic structures considered in this paper to topology will be
explored.

\begin{ackn}
We would like to thank Dennis Sullivan and Jim Stasheff for their
useful comments.
\end{ackn}

\section{Actions of chord diagrams for strict algebras}\label{strict-chapter}
This paper mainly treats chain level operations on the cyclic
Hochschild complex and the Hochschild complex of
V$_\infty$-algebras. It is useful, though, to start by describing
the operations in the simple case of an ungraded and strictly
associative algebras. We will first recall from \cite{TZ}, how
chord diagrams with marked points act on the Hochschild complex.
We then describe a cyclic version in which chord diagrams without
marked points act on the cyclic Hochschild complex of a strictly
associative algebra. We finish the section by incorporating the
two-sided cobar complex into the picture.

\subsection{Action of chord diagrams on Hochschild complex}\label{strict-Hoch-section}
\begin{figure}
\[
\begin{pspicture}(0,0.5)(7,5)
 \pscircle(1,3){0.7} \pscircle(6,4){0.7} \pscircle(4,1){0.7}
 \pscurve[linestyle=dashed](1,2.3)(1.5,1.5)(3.3,1)
 \psline[linestyle=dashed](1.6,2.6)(3.3,2.4)
 \psline[linestyle=dashed](4,1.7)(3.3,2.4)
 \psline[linestyle=dashed](5.6,3.4)(3.3,2.4)
 \pscurve[linestyle=dashed](2.6,3.4)(2,3.8)(1.4,3.6)
 \pscurve[linestyle=dashed](2.6,3.4)(2.8,1.6)(3.4,1.4)
 \pscurve[linestyle=dashed](2.6,3.4)(3.8,4.4)(5.3,4.2)
 \psline(6,4.45)(6,4.7)      \psline(4.75,0.45)(4.4,0.75)
 \psline(1.4,3.6)(1.2,3.4)  \psline(5.4,3.7)(5.1,3.6)
 \rput(1.1,3.3){\tiny$a^*$}     \rput(1.4,3.25){\tiny$a^*$}
 \rput(1.5,3){\tiny$a^*$}     \rput(1.4,2.75){\tiny$a^*$}
 \rput(1,2.52){\tiny$a^*$}    \rput(1.25,2.55){\tiny$a^*$}
 \rput(0.6,2.75){\tiny$a^*$}  \rput(0.75,2.58){\tiny$a^*$}
 \rput(0.5,3){\tiny$a^*$}     \rput(0.55,3.25){\tiny$a^*$}
 \rput(0.75,3.45){\tiny$a^*$} \rput(1,3.5){\tiny$a^*$}
 \rput(6.23,4.43){\tiny$a^*$} \rput(6.4,4.25){\tiny$a^*$}
 \rput(6.5,4){\tiny$a^*$}     \rput(6.4,3.75){\tiny$a^*$}
 \rput(6,3.52){\tiny$a^*$}    \rput(6.25,3.55){\tiny$a^*$}
 \rput(5.6,3.75){\tiny$a^*$}  \rput(5.75,3.58){\tiny$a^*$}
 \rput(5.5,4){\tiny$a^*$}     \rput(5.55,4.25){\tiny$a^*$}
 \rput(5.75,4.45){\tiny$a^*$} \rput(5.95,4.3){\tiny$a^*$}
 \rput(4.23,1.43){\tiny$a^*$} \rput(4.4,1.25){\tiny$a^*$}
 \rput(4.5,1){\tiny$a^*$}     \rput(4.3,0.76){\tiny$a^*$}
 \rput(4,0.52){\tiny$a^*$}    \rput(4.25,0.55){\tiny$a^*$}
 \rput(3.6,0.75){\tiny$a^*$}  \rput(3.75,0.58){\tiny$a^*$}
 \rput(3.5,1){\tiny$a^*$}     \rput(3.55,1.25){\tiny$a^*$}
 \rput(3.75,1.45){\tiny$a^*$} \rput(4,1.5){\tiny$a^*$}
\end{pspicture}
\]
\caption{Chord diagrams with marked points}
\label{marked-strict-graph}
\end{figure}
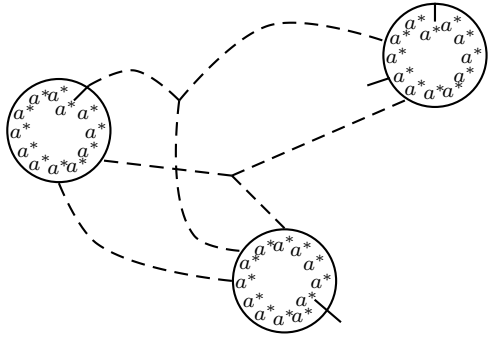

Recall from \cite{TZ} that a cyclic Sullivan chord diagram
consists of circles and chords, where the endpoints of the chords
lie on the circles. A thickening of the chord diagram will give
rise to a surface with two types of boundaries. Those situated
inside of the circles are referred to as the inputs while the
remaining ones form the outputs. Both inputs and outputs are
numerated, and there is exactly one marked point for each input
and output. It may also happen that an input marked point and an
output marked point are attached at the same spot, or that a
marked point may be attached to an endpoint of a chord. For
further details on how to treat these cases see \cite{TZ}. An
example of such a chord diagram is given in Figure
\ref{marked-strict-graph}.

Now, let $A$ be an associative algebra with unit $\1$ and an
invariant and symmetric co-inner product $U = \sum a_i\otimes
b_i\in A \otimes A $. By definition, symmetric means $\sum
a_i\otimes b_i =\sum b_i \otimes a_i$ and invariant means $\sum
a_ia \otimes b_i =\sum a_i\otimes a b_i$ and $\sum aa_i \otimes
b_i =\sum a_i\otimes b_ia$, for all $a \in A$. Note that an
invariant and symmetric inner product that is nondegenerate gives
rise to a co-inner product. The Hochschild cochain complex of $A$
with values in its dual $A^*$ is by definition the graded vector
space,
$${CH^\bullet}(A, A^*)=\prod_{j=0}^\infty \Big \{f:A^{ \otimes
j}\to A^* \Big\},$$ with the differential,
\begin{eqnarray*}
\delta f(a_1,\cdots,a_{j}) &=&  a_1 \cdot f(a_2, \cdots, a_j)\\
&&+ \sum_{i=1}^{j-2}(-1)^if(a_1, \cdots , a_ia_{i+1}, \cdots ,a_j)
\\&&+ (-1)^{j-1} f(a_1, \cdots ,a_{j-1})\cdot a_j.
\end{eqnarray*}

In the presence of a unit, one may furthermore define the
normalized subcomplex,
\begin{multline*} \overline{CH^\bullet}(A,A^*)=\prod_{j=0}^\infty
\Big\{f:A^{ \otimes j}\to A^* \Big|
\quad\quad\quad\quad\quad\quad\quad\quad\quad
\quad\quad\quad\quad\,\, \\%
f(a_1, \cdots, a_j)=0, \text{ if any of the } a_1,\cdots, a_j \text{
equals }\1 \Big\}.
\end{multline*}
The inclusion $\overline{CH^\bullet}(A,A^*) \hookrightarrow
{CH^\bullet}(A,A^*)$ is a quasi-isomorphism; see \cite{L}.

It was shown in \cite{TZ} that the PROP of cyclic Sullivan chord
diagrams with marked points acts on the normalized Hochschild
cochain complex by applying elements of
$\overline{CH^\bullet}(A,A^*)$ to the inputs of a chord diagram,
where the special last element is lined up with the input marked
point. The outputs of the chord diagram yield the output of the
operation, where the element aligned with the output marked point
is interpreted as the special last element. The role of a chord is
to apply the co-inner product to the elements at its endpoints.

\subsection{Action of chord diagrams on cyclic Hochschild complex}\label{strict-Cyc-section}
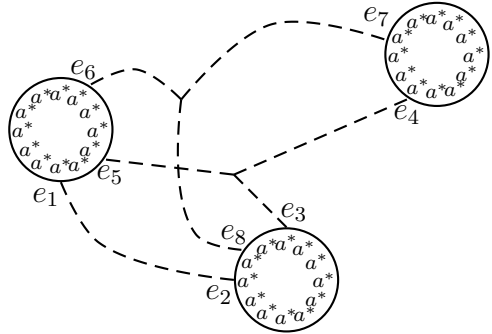
\begin{figure}
\[
\begin{pspicture}(0,0.5)(7,5)
%
 \rput(0.8,2.1){$e_1$}   \rput(1.65,2.4){$e_5$}
 \rput(3.1,0.8){$e_2$}   \rput(1.3,3.8){$e_6$}
 \rput(4.1,1.9){$e_3$}   \rput(3.3,1.6){$e_8$}
 \rput(5.6,3.2){$e_4$}   \rput(5.2,4.5){$e_7$}
 \pscircle(1,3){0.7} \pscircle(6,4){0.7} \pscircle(4,1){0.7}
 \pscurve[linestyle=dashed](1,2.3)(1.5,1.5)(3.3,1)
 \psline[linestyle=dashed](1.6,2.6)(3.3,2.4)
 \psline[linestyle=dashed](4,1.7)(3.3,2.4)
 \psline[linestyle=dashed](5.6,3.4)(3.3,2.4)
 \pscurve[linestyle=dashed](2.6,3.4)(2,3.8)(1.4,3.6)
 \pscurve[linestyle=dashed](2.6,3.4)(2.8,1.6)(3.4,1.4)
 \pscurve[linestyle=dashed](2.6,3.4)(3.8,4.4)(5.3,4.2)
 \rput(1.23,3.43){\tiny$a^*$} \rput(1.4,3.25){\tiny$a^*$}
 \rput(1.5,3){\tiny$a^*$}     \rput(1.4,2.75){\tiny$a^*$}
 \rput(1,2.52){\tiny$a^*$}    \rput(1.25,2.55){\tiny$a^*$}
 \rput(0.6,2.75){\tiny$a^*$}  \rput(0.75,2.58){\tiny$a^*$}
 \rput(0.5,3){\tiny$a^*$}     \rput(0.55,3.25){\tiny$a^*$}
 \rput(0.75,3.45){\tiny$a^*$} \rput(1,3.5){\tiny$a^*$}
 \rput(6.23,4.43){\tiny$a^*$} \rput(6.4,4.25){\tiny$a^*$}
 \rput(6.5,4){\tiny$a^*$}     \rput(6.4,3.75){\tiny$a^*$}
 \rput(6,3.52){\tiny$a^*$}    \rput(6.25,3.55){\tiny$a^*$}
 \rput(5.6,3.75){\tiny$a^*$}  \rput(5.75,3.58){\tiny$a^*$}
 \rput(5.5,4){\tiny$a^*$}     \rput(5.55,4.25){\tiny$a^*$}
 \rput(5.75,4.45){\tiny$a^*$} \rput(6,4.5){\tiny$a^*$}
 \rput(4.23,1.43){\tiny$a^*$} \rput(4.4,1.25){\tiny$a^*$}
 \rput(4.5,1){\tiny$a^*$}     \rput(4.4,0.75){\tiny$a^*$}
 \rput(4,0.52){\tiny$a^*$}    \rput(4.25,0.55){\tiny$a^*$}
 \rput(3.6,0.75){\tiny$a^*$}  \rput(3.75,0.58){\tiny$a^*$}
 \rput(3.5,1){\tiny$a^*$}     \rput(3.55,1.25){\tiny$a^*$}
 \rput(3.75,1.45){\tiny$a^*$} \rput(4,1.5){\tiny$a^*$}
\end{pspicture}
\]
\caption{Chord diagrams without marked points} \label{chord graph}
\end{figure}

Next, we will describe a cyclic version of the above action. For
this, we will look at a simpler class of chord diagrams, namely
chord diagrams similar to the ones in Figure
\ref{marked-strict-graph}, with the distinction of not having
marked points. Thus, the circles and chords determine the inputs
and outputs, but there are no longer marked points to determine
the position of a special element. An example of such a chord
diagram is shown in Figure \ref{chord graph}.
In addition to the combinatorics given by the chord diagram, we
require that each diagram comes with an orientation on the vector
space generated by the chords endpoints. For example, the chord
diagram in Figure \ref{chord graph} has eight chord endpoints,
labelled by $e_1,\cdots, e_8$, and a possible orientation could be
given by $e_1 \wedge\cdots \wedge e_8$. A chord diagram with an
orientation is identified with the negative of the same diagram
with the opposite orientation.

The boundary of a chord diagram is given by a sum obtained by
collapsing each of the circle pieces in between the endpoints of
the chords, one at a time. The induced orientation for each
collapse is determined by bringing the two chord endpoints of the
collapsing circle to the beginning of the expression for the
orientation, in the order determined by the clockwise direction of
the circle piece, and then fusing them into one point. For
example, if we would like to collapse the the circle piece between
$e_2$ and $e_8$ in Figure \ref{chord graph}, then we write $e_1
\wedge\cdots \wedge e_8=- e_2 \wedge e_8\wedge e_1\wedge
e_3\wedge\cdots \wedge e_7$ and take the induced orientation
$-e_9\wedge e_1\wedge e_3\wedge\cdots \wedge e_7$, where $e_9$ is
the chord endpoint obtained by fusing $e_2$ and $e_8$.

We may also compose two chord diagrams. This consists of
identifying the outputs of the first diagram with the inputs of
the second diagram, and then attaching the chords of the second
diagram to the first diagram in all possible ways. If
$e_1\wedge\cdots \wedge e_n$ is the orientation of the first
diagram and $e_{n+1}\wedge\cdots \wedge e_m$ is the orientation of
the second diagram, then $e_1 \wedge \cdots \wedge e_n \wedge
e_{n+1}\wedge\cdots \wedge e_m$ is the orientation of the
composition of the diagrams.

Now, let $A$ be an associative algebra concentrated in degree
zero, with an invariant and symmetric co-inner product $U$. That
is to say, $U =\sum a_i\otimes b_i \in A \otimes A$ satisfies
\begin{enumerate}
\item[$(1)$] $\sum a_i \otimes b_i =\sum b_i \otimes a_i$,
\item[$(2)$] $\sum a_ia\otimes b_i=\sum a_i\otimes ab_i$, for
every $a \in A$, \item[$(3)$] $\sum aa_i \otimes b_i=\sum
a_i\otimes b_ia$, for every $a\in A$.
\end{enumerate}
We consider the cyclic Hochschild cochain complex,
\begin{multline*}
CC^\bullet(A)=\prod_{j=0}^\infty \Big\{f:A^{\otimes j+1}\to {\bf k}
\Big|
\\ f(a_1,\cdots,a_{j+1})=(-1)^{j} \cdot f(a_{j+1},a_1,\cdots, a_j)
\Big\},
\end{multline*}
with the differential,
\begin{eqnarray*}
\delta f(a_0, \cdots, a_j) &=& \sum_{i=0}^{j-1}  (-1)^i f(a_0,
\cdots, a_ia_{i+1}, \cdots, a_j) \\&& +(-1)^j (a_ja_0, \cdots,
a_{j-1}).
\end{eqnarray*}
To better understand the signs in the above expression, it would
be more conceptual to think that elements of $A$ have been shifted
to degree $-1$, and that the signs are those coming from by the
Koszul sign rule which introduces a $(-1)^{pq}$, whenever
something of degree $p$ is moved passed something of degree $q$.

We claim that the PROP of chord diagrams without marked points acts
on $CC^\bullet(A)$. In fact, this action is very similar to the
action on $\overline{CH^\bullet}(A,A^*)$ from the previous section.
Elements $f_1, \cdots, f_k$ of $CC^\bullet(A)$ will be applied to
the input circles of the chord diagram. We may write the inputs
$f_1\otimes \cdots\otimes f_k$ in a linear order staring from the
first to the $k^{th}$ input. Then, much like before, for each chord,
the co-inner product is applied to the elements at the chords
endpoints. To obtain the correct signs, before applying the co-inner
product to the elements of the shifted $A^*$, one needs to shift
them back to degree zero. In order to unshift them, we apply shift
operators of degree $1$ to each chord endpoint. More precisely, we
apply shift operations on the left of $f_1\otimes \cdots\otimes f_k$
in the order given by the orientation $e_1\wedge \cdots\wedge e_n$
of the diagram. Then, the inner products for the chords may safely
be applied, since the chord endpoints will now have degree zero. We
collect the outputs of the chord diagram as outputs of the
operation. Since these outputs will in general not be cyclically
invariant, we perform a final cyclic symmetrization to each output.

\subsection{Action of open/closed chord diagrams on cyclic
Hochschild complex coupled with two-sided cobar complex}
\label{strict-open-closed-section} 
\begin{figure}
\[
\begin{pspicture}(-2,0.5)(9,5)
 \pscircle(1,3){0.7} \pscircle(6,4){0.7} \pscircle(4,1){0.7}
 \pscurve[linestyle=dashed](1,2.3)(1.5,1.5)(3.3,1)
 \psline[linestyle=dashed](1.6,2.6)(3.3,2.4)
 \psline[linestyle=dashed](4,1.7)(3.3,2.4)
 \psline[linestyle=dashed](5.6,3.4)(3.3,2.4)
 \pscurve[linestyle=dashed](2.6,3.4)(2,3.8)(1.4,3.6)
 \pscurve[linestyle=dashed](2.6,3.4)(2.8,1.6)(3.4,1.4)
 \pscurve[linestyle=dashed](2.6,3.4)(3.8,4.4)(5.3,4.2)
 \psline(-1.2,2.3)(-1.2,4.3) \psline(6.6,1)(7,2.8)
 \psline[linestyle=dashed](-1.2,3)(0.3,3)
 \pscurve[linestyle=dashed](-1.2,3.6)(-0.2,4)(1,4.2)(1.8,4)(1.7,3.2)
 \pscurve[linestyle=dashed](6.8,1.9)(5.5,2.3)(3.3,2.4)
 \rput(1.23,3.43){\tiny$a^*$} \rput(1.4,3.25){\tiny$a^*$}
 \rput(1.5,3){\tiny$a^*$}     \rput(1.4,2.75){\tiny$a^*$}
 \rput(1,2.52){\tiny$a^*$}    \rput(1.25,2.55){\tiny$a^*$}
 \rput(0.6,2.75){\tiny$a^*$}  \rput(0.75,2.58){\tiny$a^*$}
 \rput(0.5,3){\tiny$a^*$}     \rput(0.55,3.25){\tiny$a^*$}
 \rput(0.75,3.45){\tiny$a^*$} \rput(1,3.5){\tiny$a^*$}
 \rput(6.23,4.43){\tiny$a^*$} \rput(6.4,4.25){\tiny$a^*$}
 \rput(6.5,4){\tiny$a^*$}     \rput(6.4,3.75){\tiny$a^*$}
 \rput(6,3.52){\tiny$a^*$}    \rput(6.25,3.55){\tiny$a^*$}
 \rput(5.6,3.75){\tiny$a^*$}  \rput(5.75,3.58){\tiny$a^*$}
 \rput(5.5,4){\tiny$a^*$}     \rput(5.55,4.25){\tiny$a^*$}
 \rput(5.75,4.45){\tiny$a^*$} \rput(6,4.5){\tiny$a^*$}
 \rput(4.23,1.43){\tiny$a^*$} \rput(4.4,1.25){\tiny$a^*$}
 \rput(4.5,1){\tiny$a^*$}     \rput(4.4,0.75){\tiny$a^*$}
 \rput(4,0.52){\tiny$a^*$}    \rput(4.25,0.55){\tiny$a^*$}
 \rput(3.6,0.75){\tiny$a^*$}  \rput(3.75,0.58){\tiny$a^*$}
 \rput(3.5,1){\tiny$a^*$}     \rput(3.55,1.25){\tiny$a^*$}
 \rput(3.75,1.45){\tiny$a^*$} \rput(4,1.5){\tiny$a^*$}
 \rput(-1.4,4.3){\tiny$m^*$}  \rput(-1.4,3.9){\tiny$a^*$}
 \rput(-1.4,3.5){\tiny$a^*$}  \rput(-1.4,3.1){\tiny$a^*$}
 \rput(-1.4,2.7){\tiny$a^*$}  \rput(-1.4,2.3){\tiny$n^*$}
 \rput(7.30,2.80){\tiny$n^*$}  \rput(7.14,2.45){\tiny$a^*$}
 \rput(7.08,2.10){\tiny$a^*$}  \rput(7.02,1.70){\tiny$a^*$}
 \rput(6.96,1.35){\tiny$a^*$}  \rput(6.90,1.00){\tiny$m^*$}
\end{pspicture}
\]
\caption{Open/closed chord diagrams without marked points}
\label{open-closed-strict-graph}
\end{figure}
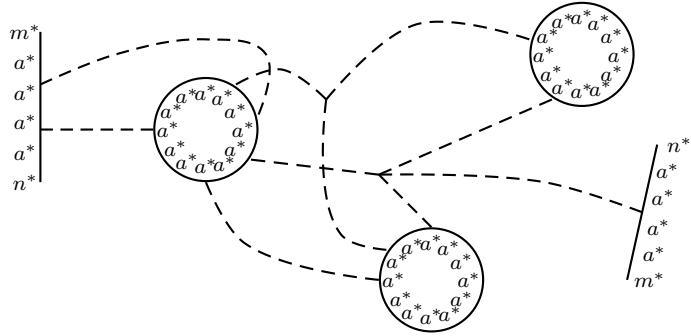

As yet another variation, we
describe an action of the open/closed chord diagrams. For this, we
look at chord diagrams without marked points, with some of the
circles opened up to intervals. Thus, the open/closed chord
diagrams consist of intervals, circles and chords; see Figure
\ref{open-closed-strict-graph}.
Since there are no marked points, we additionally require an
orientation on the vector space generated by the chords endpoints.

Once again, let $A$ be an associative algebra with an invariant
and symmetric co-inner product $U\in A \otimes A$, and let $M$ and
$N$ be left and right modules over $A$, respectively. Consider the
two-sided cobar complex,
$$ C^\bullet(M,A,N)= \prod_{j=0}^\infty \Big\{f:M\otimes
A^{\otimes j}\otimes N\to {\bf k}\Big\}$$ Using the left and right
module structures on $M$ and $N$, one defines a differential on
$C^\bullet(M,A,N)$; see section \ref{open section}.

We claim that the PROP of open/closed chord diagrams acts on
$CC^\bullet(A)\otimes C^\bullet(M,A,N)$. Much like before,
$CC^\bullet(A)$ corresponds to the closed circles, and
$C^\bullet(M,A,N)$ to the open intervals. To be more precise, we
associate $M$ to the left endpoint and $N$ to the right endpoint
of the interval. Once again, chords will apply co-inner products,
and we can collect the outputs according to the outputs of the
diagram. An orientation on the vector space generated by the
chords endpoints is needed to determine the order in which the
inverse of the shift operator, bringing elements of $A^*$ from the
shifted degree $+1$ back to degree $0$, must be applied.

\section{V$_\infty$-algebras}\label{algebra section}

One would like to extend the above discussions to a category in
which associative algebras are replaced by A$_\infty$-algebras.
Although the notion of homotopy co-inner products proves helpful
in the study of Hochschild cohomology $HH^\bullet(A, A^*)$ of an
A$_\infty$-algebra $A$ (see \cite{T2}) there are limitations when
it comes to Hochschild cochains $CH^\bullet(A, A^*)$. These
limitations seem to be more pronounced in the cyclic case of
$CC^\bullet(A)$, as we will point out in section
\ref{obstract-and-operat}. In fact, in order to see that the
cyclic Hochschild cohomology is a Lie algebra, higher homotopies
that are not part of the homotopy co-inner product are needed to
guarantee the Jacobi identity.

In this section, we describe the concept of a V$_k$-algebra, which
is germane to the study of the algebraic structure of the
Hochschild complex and cyclic Hochschild complex of an
A$_\infty$-algebra. The notion of an A$_\infty$-algebra is only
the first layer of this structure. The second layer is where the
concept of an invariant and symmetric homotopy co-inner product
comes into its own. The third and higher layers are new concepts
that we shall develop and use. We will see in section \ref{action
section} how considering V$_k$-algebras will shed some light on
the algebraic structure of the cyclic Hochschild complex and its
cousins.

\subsection{A$_\infty$-algebras or $V_1$-algebras}\label{A-infty-section}

\begin{figure}
\[
     {\psset{arrows=<-}
    \pstree[treemode=U, levelsep=0.8cm]{\Tp}
    {\pstree{\Tc*{3pt}}{\pstree{\Tc*{3pt}} {\Tp \Tp}\pstree{\Tn}{\Tp}}
    }}
\begin{pspicture}(0,0)(1,2) \rput(.5,1.3){$=$} \end{pspicture}
     \psset{arrows=<-}
     \pstree[treemode=U, levelsep=0.8cm]{\Tp}
    {\pstree{\Tc*{3pt}}{\pstree{\Tn}{\Tp}\pstree{\Tc*{3pt}} {\Tp \Tp}}}
\] \caption{Associativity}\label{associativity}
\end{figure}
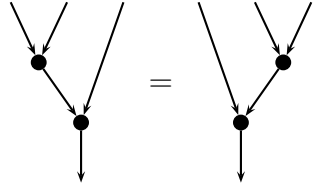

Let $\left(A=\bigoplus_{n\in \mathbb Z}A_n,\partial_A\right)$ be a
complex over a field ${\bf k}$ of characteristic zero, such that
each $A_n$ is a finite dimensional vector space. Consider the dual
complex $\left(A^\ast= \bigoplus_{n\in \mathbb
Z}A^*_n,\partial_A^*\right)$. Note that the grading for the dual
space $A^\ast$ is negative of that for $A$. For instance, if $A$ is
concentrated in negative degrees, then $A^*$ is positively graded.
The finite dimensionality assumption allows for such identifications
as $Hom(A, A)= A^* \otimes A$ as graded objects. An associative
algebra structure on $A$ consists of a chain map $\mu_2:A^{\otimes
2}\to A$, usually written as $\mu_2(a_1,a_2)=a_1\cdot a_2$, such
that $(a_1\cdot a_2)\cdot a_3= a_1\cdot (a_2\cdot a_3)$. A
representation of the identities for associativity using graphs is
shown in Figure \ref{associativity}, where incoming edges denote the
arguments of $\mu_2$ and an outgoing edge stands for its output.

\begin{figure}
\[ \begin{pspicture}(0,0)(1,2) \rput(.5,1){$\sum$} \end{pspicture}
\psset{arrows=<-} \pstree[treemode=U, levelsep=0.8cm,
treesep=0.3cm]{\Tp} {\pstree{\Tc*{3pt}}{\Tp \Tp\Tp\pstree{\Tc*{3pt}}
{\Tp\Tp\Tp\Tp \Tp}\Tp\Tp} }
\begin{pspicture}(0,0)(1,2) \rput(.5,1){$=0$} \end{pspicture}
\] \caption{Homotopy associativity (A$_\infty$)}\label{a-infty}
\end{figure}
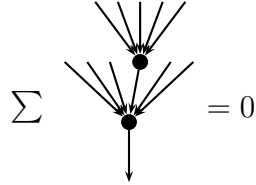

The concept of an A$_\infty$-algebra, introduced by J. Stasheff in
\cite{St}, is a generalization of the concept of associative
algebra. An A$_\infty$-algebra is a vector space $A$ endowed with
a sequence of maps $\{\mu_n:A^{ \otimes n}\to A\}_{n \geq 1}$,
where the sum of all compositions of any $\mu_i$ and $\mu_j$ is
zero, $$ \sum_{i,j,k} \pm \mu_i(a_1, \cdots, \mu_j(a_k,
\cdots,a_{k+j-1}),\cdots,a_n)=0. $$ We defer the discussion of
signs until section \ref{vert O, I, II section}, where we cover
the general case. A graphical representation of the above
identities can be seen in Figure \ref{a-infty}.
Here, the $j$ incoming edges stand for the arguments of $\mu_j$
and the one outgoing edge is used for its output. Put slightly
differently, we may view $\mu_j:A^{\otimes j}\to A$ as an element
$\mu_j\in A \otimes (A^*)^{\otimes j}$, and reserve the incoming
edges for the elements of $A^*$, and the outgoing edge for those
of $A$. This assignment of elements to the edges, which is
suitable for evaluating of elements of $A^*$ on $A$, is used
throughout the paper.

\subsection{A$_\infty$-algebras with homotopy co-inner products or
$V_2$-algebras}\label{homotopy co-inner product-section}

We also need a generalization of the concept of an invariant and
symmetric co-inner product. By definition, an invariant and
symmetric co-inner product of dimension $d$ on a graded associative
algebra $A$ (finite dimensional in each degree, and without a unit)
is an element $U =\sum a_i\otimes b_i \in A \otimes A$ of degree
$-d$, called the co-inner product, satisfying,
\begin{enumerate}
\item[$(1)$] $\sum a_i \otimes b_i =\sum \pm b_i \otimes a_i$,
\item[$(2)$] $\sum a_ia\otimes b_i=\sum \pm a_i\otimes ab_i$, for
every $a \in A$, \item[$(3)$] $\sum aa_i \otimes b_i=\sum \pm
a_i\otimes b_ia$, for every $a\in A$,
\end{enumerate}
where the signs $\pm$ are those dictated by the Koszul sign rule.

In practice, co-inner products of dimension $d$ are often obtained
from inner products $\langle\cdot, \cdot \rangle:A\otimes A \to {\bf
k}$, of degree $d$, such that,
\begin{itemize}
\item[$(1')$] $\langle a, b\rangle=\langle a, b\rangle$,
\item[$(2')$] $\langle ab,c\rangle=\langle a, bc\rangle$,
\item[$(3')$] $\langle a, bc\rangle =\langle ca, b\rangle$,
\item[$(4')$] $A \to A^\ast$, defined by $a \mapsto \langle a,
\cdot \rangle$,  is an isomorphism, for every $a\in A$.
\end{itemize}
Thus, sing the isomorphism $A \cong A^\ast$ one transports the
inner product $\langle \cdot, \cdot\rangle \in A^\ast \otimes
A^\ast$ to get an element of $U \in A \otimes A$.

For instance, for a Poincar\'e duality space $X$ of dimension $d$,
consider the cohomology $H^\bullet X$, negatively graded, and
define $\langle \cdot , \cdot \rangle: H^\bullet X \times
H^\bullet X \rightarrow {\bf k}$ by $\langle a, b \rangle = (a
\cup b)[X]$. The co-inner product is given by the Thom class of
the diagonal $U \in H^\bullet X \otimes H^\bullet X$.

\begin{figure}
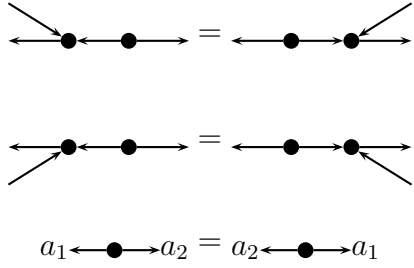

\[
    \psset{arrows=<-}\pstree[treemode=L, levelsep=0.8cm, treesep=0.5cm]{\Tp}
    {\pstree[arrows=->]{\Tc*{3pt}}{ \pstree{\Tc*{3pt}}
    {\pstree[levelsep=0, arrows=<-]{\Tn}{\Tp}\Tp\Tn}} }
  =
    \psset{arrows=<-}\pstree[treemode=R, levelsep=0.8cm, treesep=0.5cm]{\Tp}
    {\pstree[arrows=->]{\Tc*{3pt}}{ \pstree{\Tc*{3pt}}
    {\pstree[levelsep=0, arrows=<-]{\Tn}{\Tp}\Tp\Tn}} }
 \]

 \[
    \psset{arrows=<-}\pstree[treemode=L, levelsep=0.8cm, treesep=0.5cm]{\Tp}
    {\pstree[arrows=->]{\Tc*{3pt}}{ \pstree{\Tc*{3pt}}
    {\Tn\Tp\pstree[levelsep=0, arrows=<-]{\Tn}{\Tp}}} }
  =
    \psset{arrows=<-}\pstree[treemode=R, levelsep=0.8cm, treesep=0.5cm]{\Tp}
    {\pstree[arrows=->]{\Tc*{3pt}}{ \pstree{\Tc*{3pt}}
    {\Tn\Tp\pstree[levelsep=0, arrows=<-]{\Tn}{\Tp}}} }
 \]

 \[
    \psset{arrows=<-}\pstree[treemode=L, levelsep=0.8cm]{\Tr{a_2}}
    {\pstree[arrows=->]{\Tc*{3pt}}{ \Tr{a_1}} }
  =
    \psset{arrows=<-}\pstree[treemode=L, levelsep=0.8cm]{\Tr{a_1}}
    {\pstree[arrows=->]{\Tc*{3pt}}{ \Tr{a_2}} }
 \]
\caption{Co-inner product}\label{co-inner product}
\end{figure}

In terms of graphs, a co-inner product $U \in A\otimes A$ is
represented by a vertex with two outgoing edges. Figure
\ref{co-inner product} represents the required identities in terms
of graphs.
Notice that the cyclic ordering of the arguments $a^*_1, a^*_2$
and $a$ are preserved in each of the equations.

Just as associative algebras were generalized to
A$_\infty$-algebras, co-inner products may be generalized to
homotopy co-inner products. A treatment of such a generalization for
any cyclic operad appears in \cite{LT}. As it turns out, the
homotopy co-inner products are given by elements $v_{i, j}\in
A\otimes (A^*)^{\otimes i}\otimes A\otimes (A^*)^{\otimes j}$. These
$v_{i, j}$'s satisfy generalized co-inner product conditions
indicated in Figure \ref{homotopy co-inner product}.
\begin{figure}
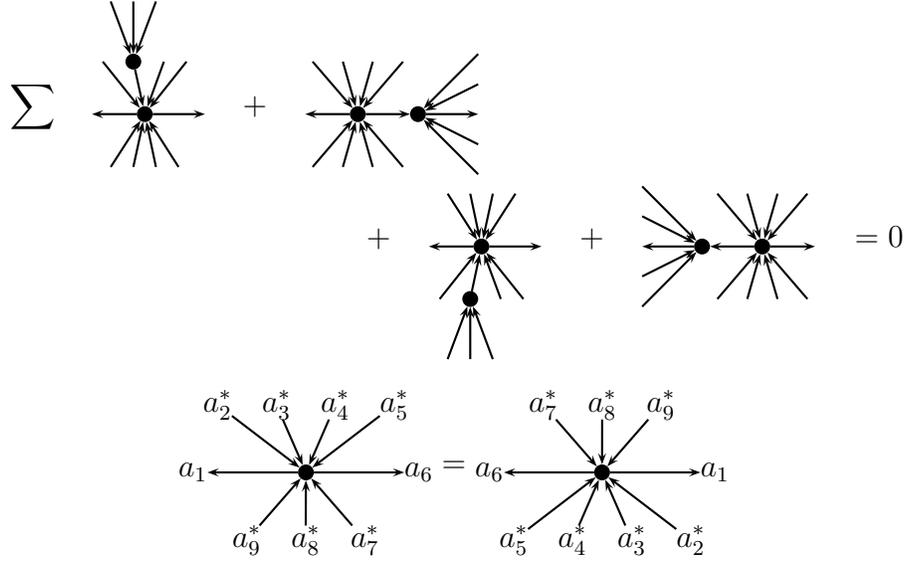

\begin{multline*}
\sum\quad
    \psset{arrows=<-}
    \pstree[treemode=R, levelsep=0.7cm, treesep=0.3cm]{\Tp}{
    \pstree[levelsep=0]{\Tc*{3pt}}{
    \pstree[levelsep=0.4cm,treemode=U]{\Tn}{\Tp
         \pstree[levelsep=0.8cm]{\Tc*{3pt}}{\Tp\Tp\Tp}\Tp\Tp}
    \pstree[levelsep=0.8cm,treemode=R, arrows=->]{\Tn}{\Tp}
    \pstree[levelsep=0.4cm,treemode=D]{\Tn}{\Tp\Tp\Tp\Tp}
    }}\quad + \quad
    \psset{arrows=<-}
    \pstree[treemode=R, levelsep=0.7cm, treesep=0.4cm]{\Tp}{
    \pstree[levelsep=0]{\Tc*{3pt}}{
    \pstree[levelsep=0.3cm,treemode=U]{\Tn}{\Tp\Tp\Tp\Tp}
    \pstree[levelsep=0.8cm,treemode=R, arrows=->]{\Tn}{
         \pstree[arrows=<-]{\Tc*{3pt}}{\Tp\Tp\pstree[arrows=->, levelsep=0]{\Tn}{\Tp}\Tp\Tp}}
    \pstree[levelsep=0.3cm,treemode=D]{\Tn}{\Tp\Tp\Tp\Tp}
    }}\\ + \quad
    \psset{arrows=<-}
    \pstree[treemode=R, levelsep=0.7cm, treesep=0.3cm]{\Tp}{
    \pstree[levelsep=0]{\Tc*{3pt}}{
    \pstree[levelsep=0.4cm,treemode=U]{\Tn}{\Tp\Tp\Tp\Tp}
    \pstree[levelsep=0.8cm,treemode=R, arrows=->]{\Tn}{\Tp}
    \pstree[levelsep=0.4cm,treemode=D]{\Tn}{\Tp
       \pstree[levelsep=0.8cm]{\Tc*{3pt}}{\Tp\Tp\Tp}\Tp\Tp}
    }} \quad + \quad
    \psset{arrows=<-}
    \pstree[treemode=L, levelsep=0.7cm, treesep=0.4cm]{\Tp}{
    \pstree[levelsep=0]{\Tc*{3pt}}{
    \pstree[levelsep=0.3cm,treemode=U]{\Tn}{\Tp\Tp\Tp\Tp}
    \pstree[levelsep=0.8cm,treemode=L, arrows=->]{\Tn}{
         \pstree[arrows=<-]{\Tc*{3pt}}{\Tp\Tp\pstree[arrows=->, levelsep=0]{\Tn}{\Tp}\Tp\Tp}}
    \pstree[levelsep=0.3cm,treemode=D]{\Tn}{\Tp\Tp\Tp\Tp}
    }}\quad= 0
\end{multline*}

\[
    \psset{arrows=<-}
    \pstree[treemode=R, levelsep=1.5cm, treesep=0.4cm]{\Tr{a_1}}{
    \pstree[levelsep=0]{\Tc*{3pt}}{
    \pstree[levelsep=0.5cm,treemode=U]{\Tn}{\Tr{a^*_2}\Tr{a^*_3}\Tr{a^*_4}\Tr{a^*_5}}
    \pstree[levelsep=1.5cm,treemode=R, arrows=->]{\Tn}{\Tr{a_6}}
    \pstree[levelsep=0.5cm,treemode=D]{\Tn}{\Tr{a^*_9}\Tr{a^*_8}\Tr{a^*_7}}
    }}
    =
    \psset{arrows=<-}
    \pstree[treemode=L, levelsep=1.5cm, treesep=0.4cm]{\Tr{a_1}}{
    \pstree[levelsep=0]{\Tc*{3pt}}{
    \pstree[levelsep=0.5cm,treemode=U]{\Tn}{\Tr{a^*_7}\Tr{a^*_8}\Tr{a^*_9}}
    \pstree[levelsep=1.5cm,treemode=L, arrows=->]{\Tn}{\Tr{a_6}}
    \pstree[levelsep=0.5cm,treemode=D]{\Tn}{\Tr{a^*_5}\Tr{a^*_4}\Tr{a^*_3}\Tr{a^*_2}}
    }}
\]
\caption{Homotopy co-inner product} \label{homotopy co-inner
product}
\end{figure}
In Figure \ref{homotopy co-inner product}, the upper equation is
the generalized bimodule condition. Here, each $v_{i, j}$ is
pictured by a vertex with two (horizontal) outgoing edges, $i$
incoming edges on top, and $j$ incoming edges at the bottom. For
each $i,j\geq 0$, we go around $v_{i, j}$ and for all $m \geq 1$
attach the maps $\mu_m$ of the A$_\infty$-algebra structure in all
possible ways. Notice the difference between attaching a $\mu_m$
to one of the horizontal arrows as opposed to attaching it to one
from the top or the bottom. It is also worth mentioning that the
above sum splits into countably many equations, each of which is a
finite sum with a fixed number of edges. The lower equation in
Figure \ref{homotopy co-inner product} shows how the symmetry
condition is generalized. It is simply by requiring that each
$v_{i, j}$ is related to $v_{j, i}$ via a $180^\circ$ rotation.

\subsection{V$_k$-algebras, for ~ $ {\bf 1\leq k \leq \infty}$}
\label{vert O, I, II section}

In previous sections, we looked at algebraic objects that were
labelled graphically by vertices with an arbitrary number of
incoming edges but only one or at most two outgoing ones. We call
these vertices of types $1$ and $2$, respectively. More generally,
a vertex of type $n$ has an unrestricted number of incoming edges
shuffled in between $n$ outgoing edges.
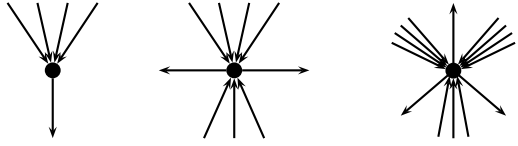
\begin{figure}
\[
    \psset{arrows=<-}
    \pstree[treemode=L, levelsep=1cm, treesep=0.4cm]{\Tn}{
    \pstree[levelsep=0]{\Tc*{3pt}}{
    \pstree[levelsep=0.7cm,treemode=U]{\Tn}{\Tp\Tp\Tp\Tp}
    \pstree[levelsep=0.7cm,treemode=D,arrows=->]{\Tn}{\Tp}
    }}
\quad
    \psset{arrows=<-}
    \pstree[treemode=L, levelsep=1cm, treesep=0.4cm]{\Tp}{
    \pstree[levelsep=0]{\Tc*{3pt}}{
    \pstree[levelsep=0.5cm,treemode=U]{\Tn}{\Tp\Tp\Tp\Tp}
    \pstree[levelsep=1cm,treemode=L,arrows=->]{\Tn}{\Tp}
    \pstree[levelsep=0.5cm,treemode=D]{\Tn}{\Tp\Tp\Tp}
    }}
\quad
\begin{pspicture}(-0.5,1)(2,2)
 \psdots[dotsize=6pt](1,1)
 \psline[arrows=->](1,1)(0.3,0.4) \psline[arrows=->](1,1)(1,1.9)
 \psline[arrows=->](1,1)(1.7,0.4) \psline[arrows=->](1,0.1)(1,0.9)
 \psline[arrows=->](0.8,0.12)(0.95,0.93) \psline[arrows=->](1.2,0.12)(1.05,0.93)
 \psline[arrows=->](1.7,1.6)(1.08,1.06)  \psline[arrows=->](1.6,1.7)(1.06,1.08)
 \psline[arrows=->](1.82,1.37)(1.1,1.02) \psline[arrows=->](1.78,1.5)(1.09,1.04)
 \psline[arrows=->](0.3,1.6)(0.92,1.06)  \psline[arrows=->](0.4,1.7)(0.94,1.08)
 \psline[arrows=->](0.18,1.37)(0.9,1.02) \psline[arrows=->](0.22,1.5)(0.91,1.04)
\end{pspicture}
\] \caption{Vertices of type $1$, $2$ and $3$} \label{type O, I,
II}
\end{figure}
The incoming edges correspond to elements of $A^*$, whereas the
ones outgoing are labelled with elements of $A$.

\begin{eqnarray*}
\text{type } 1 \leftrightarrow &\mu_j \in A \otimes (A^*)^{\otimes
j} & \leftrightarrow A_\infty \text{-algebra}\\
\text{type } 2 \leftrightarrow &v_{i, j}\in A \otimes
(A^*)^{\otimes i}\otimes A \otimes (A^* )^{ \otimes j}
&\leftrightarrow \text{Homotopy}\\&& \quad\quad \text{co-inner
product}
\end{eqnarray*}
We now consider algebraic objects that correspond to vertices of
type $n$, for $n \in \mathbb{N}$. Each vertex of type $n$ labels
an element $v_{i_1,\cdots,i_n}\in A \otimes (A^*)^{\otimes i_1}
\otimes \cdots \otimes A \otimes (A^*)^{\otimes i_n}$ for
$i_1,\cdots, i_n \in \mathbb{N}\cup \{0\}$, satisfying conditions
analogous to those in Figure \ref{homotopy co-inner product}.
\begin{defn}\label{V_k}
For fixed $k\in \{1, 2,3,\cdots, \infty\}$, a V$_k$-algebra
structure of dimension $d$ on $A$ consists of a sequence of
elements,
\begin{eqnarray*}
v_{i_1} &\in& A \otimes (A^*)^{\otimes i_1}  \quad , i_1 \neq 0\\ %
v_{i_1,i_2} &\in& A \otimes (A^*)^{\otimes i_1}\otimes A\otimes
(A^*)^{\otimes
i_2}\\%
&\vdots & \\%
v_{i_1,\cdots,i_k}&\in& A \otimes (A^*)^{\otimes i_1}
\otimes\cdots \otimes A \otimes (A^*)^{\otimes i_k}
\end{eqnarray*}
for any indices $i_1,\cdots, i_n\in \mathbb{N} \cup \{0\}$, where
$1\leq n < k+1$, subject to the following conditions.
\begin{itemize}

\item \emph{Degree:} Each $v_{i_1,\cdots,i_k}$ is an element of
degree $n(2-d)+(d-4)$, with respect to the grading of $A\otimes
(A^*)^{\otimes i_1} \otimes\cdots \otimes A \otimes (A^*)^{\otimes
i_k}$ obtained by adding $-1$ to the degrees of the elements of
$A^*$, but leaving the degree of elements of $A$ unshifted.

\item \emph{Symmetry condition:} If
$\tau:A\otimes (A^*)^{\otimes i_1}\otimes A \otimes(A^*)^{\otimes
i_2}\otimes A \otimes \cdots \otimes A \otimes (A^*)^{\otimes
i_n}\to A\otimes (A^*)^{ \otimes i_n}\otimes A
\otimes(A^*)^{\otimes i_1}\otimes \cdots \otimes A \otimes
(A^*)^{\otimes i_{n-1}}$ denotes the cyclic rotation of tensor
factors, then,
$$ v_{i_2,\cdots,i_n,i_1}=(-1)^\epsilon\cdot \tau(v_{i_1,i_2,
\cdots,i_n}) $$ Here, $(-1)^\epsilon$ is obtained by the Koszul
sign rule in the above grading; see Figure
\ref{symmetry-condition}.
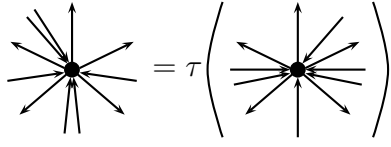
\begin{figure}
\[
\begin{pspicture}(0,0)(2,2)
 \psdots[dotsize=6pt](1,1)
 \psline[arrows=->](1,1)(0.3,0.4) \psline[arrows=->](1,1)(1,1.9)
 \psline[arrows=->](1,1)(1.7,0.4)
 \psline[arrows=<-](1.82,1.37)(1.1,1.02)  \psline[arrows=<-](0.18,1.37)(0.9,1.02)
 \psline[arrows=->](0.4,1.7)(0.94,1.08) \psline[arrows=->](0.5,1.8)(0.96,1.09)
 \psline[arrows=->](0.9,0.15)(0.98,0.9) \psline[arrows=->](1.1,0.15)(1.02,0.9)
 \psline[arrows=->](0.14,0.85)(0.9,0.96) \psline[arrows=->](1.86,0.85)(1.1,0.96)
\end{pspicture}
\begin{pspicture}(-1,0)(2.4,2)
 \psdots[dotsize=6pt](1,1) \rput(-0.6,1){$=\tau$}
 \psline[arrows=->](1,1)(0.3,0.4) \psline[arrows=->](1,1)(1,1.9)
 \psline[arrows=->](1,1)(1.7,0.4) \psline[arrows=->](1,0.1)(1,0.9)
 \psline[arrows=->](1.6,1.7)(1.06,1.08) \psline[arrows=<-](1.82,1.37)(1.1,1.02)
 \psline[arrows=<-](0.18,1.37)(0.9,1.02)
 \psline[arrows=->](0.1,1)(0.9,1)  \psline[arrows=->](1.9,1)(1.1,1)
 \psline[arrows=->](0.15,0.8)(0.9,0.95)  \psline[arrows=->](1.85,0.8)(1.1,0.95)
 \pscurve(0,1.9)(-0.2,1)(0,0.1)  \pscurve(2,1.9)(2.2,1)(2,0.1)
\end{pspicture}
\] \caption{The symmetry condition}
 \label{symmetry-condition}
\end{figure}

\item \emph{Boundary condition:} For fixed $1\leq n < k+1$, and
$i_1,\cdots, i_n \in \mathbb{N} \cup \{0\}$, let
$\Gamma=\Gamma(i_1,\cdots,i_n)$ be the set of rooted and directed
trees, that are obtained by expanding exactly one edge in the
unique tree given by the type $n$ vertex with $m=i_1+ \cdots +i_n$
incoming edges; see Figure \ref{one-internal-edge}.
\begin{figure}
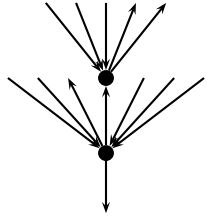

\[
    \pstree[treemode=U, levelsep=0.8cm, treesep=0.4cm, arrows=<-]{\Tp}{
    \pstree[levelsep=0.5cm,arrows=->]{\Tc*{3pt}}{
    \pstree[levelsep=0.5cm,arrows=<-]{\Tn}{\Tp }
    \pstree[levelsep=0.5cm,arrows=<-]{\Tn}{\Tp }
    \pstree[levelsep=0.5cm,arrows=->]{\Tn}{\Tp }
    \pstree[levelsep=0.5cm,arrows=->]{\Tn}{
        \pstree[levelsep=0.5cm,arrows=->]{\Tc*{3pt}}{
    \pstree[levelsep=0.5cm,arrows=<-]{\Tn}{\Tp }
    \pstree[levelsep=0.5cm,arrows=<-]{\Tn}{\Tp }
    \pstree[levelsep=0.5cm,arrows=<-]{\Tn}{\Tp }
    \pstree[levelsep=0.5cm,arrows=->]{\Tn}{\Tp }
    \pstree[levelsep=0.5cm,arrows=->]{\Tn}{\Tp }
    } }
    \pstree[levelsep=0.5cm,arrows=<-]{\Tn}{\Tp }
    \pstree[levelsep=0.5cm,arrows=<-]{\Tn}{\Tp }
    \pstree[levelsep=0.5cm,arrows=<-]{\Tn}{\Tp }
    }}
\] \caption{A tree in $\Gamma=\Gamma(2,3,0,3)$, with $n=4$ outgoing edges}
\label{one-internal-edge}
\end{figure}
For a tree $\gamma\in\Gamma$, there is no restriction on the
orientation of the internal edge (compare for example the first two
terms in Figure \ref{boundary-condition}). Also, $\gamma$ has
exactly two internal vertices, with $j_1,\cdots,j_r$ and
$j'_1,\cdots,j'_s$ incoming edges, respectively. We can assign
$v_{j_1,\cdots,j_r}$ and $v_{j'_1,\cdots,j'_s}$ to these vertices,
and write them in the order given by the internal edge of $\gamma$,
say $v_{j_1,\cdots,j_r}\otimes v_{j'_1,\cdots,j'_s}$. Then, we fuse
one tensor factor of $v_{j_1,\cdots,j_r}$ with one tensor factor of
$v_{j'_1,\cdots,j'_s}$, as indicated by the internal edge of
$\gamma$. The result is denoted by $(-1)^\epsilon \cdot v(\gamma)\in
A \otimes (A^*)^{\otimes i_1} \otimes\cdots \otimes A \otimes
(A^*)^{\otimes i_n}$, where $\epsilon$ is obtained by the product of
shifted degrees that have switched positions in the linear expansion
of $v_{j_1,\cdots,j_r}\otimes v_{j'_1,\cdots,j'_s}$. In this
notation, we require the following boundary condition; see Figure
\ref{boundary-condition}.
$$ \sum_{\gamma\in \Gamma} (-1)^\epsilon\cdot v(\gamma)=0 $$
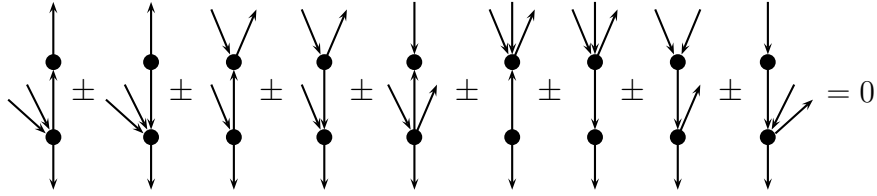
\begin{figure}
\[
\begin{pspicture}(-0.2,0)(1.1,2.5)
 \psdots[dotsize=6pt](0.5,0.8)  \psdots[dotsize=6pt](0.5,1.8)
 \psline[arrows=->](0.5,0.8)(0.5,0.1) \rput(0.9,1.4){$\pm$}
 \psline[arrows=<-](0.5,1.7)(0.5,0.8) \psline[arrows=->](0.5,1.8)(0.5,2.6)
 \psline[arrows=->](0.15,1.5)(0.45,0.9) \psline[arrows=->](-0.1,1.3)(0.4,0.85)
\end{pspicture}
\begin{pspicture}(-0.2,0)(1.1,2.5)
 \psdots[dotsize=6pt](0.5,0.8)  \psdots[dotsize=6pt](0.5,1.8)
 \psline[arrows=->](0.5,0.8)(0.5,0.1) \rput(0.9,1.4){$\pm$}
 \psline[arrows=->](0.5,1.8)(0.5,0.9) \psline[arrows=->](0.5,1.8)(0.5,2.6)
 \psline[arrows=->](0.15,1.5)(0.45,0.9) \psline[arrows=->](-0.1,1.3)(0.4,0.85)
\end{pspicture}
\begin{pspicture}(0.0,0)(1.2,2.5)
 \psdots[dotsize=6pt](0.5,0.8)  \psdots[dotsize=6pt](0.5,1.8)
 \psline[arrows=->](0.5,0.8)(0.5,0.1) \rput(1.0,1.4){$\pm$}
 \psline[arrows=->](0.5,0.8)(0.5,1.7) \psline[arrows=->](0.2,1.5)(0.45,0.9)
 \psline[arrows=->](0.2,2.5)(0.45,1.9) \psline[arrows=->](0.5,1.8)(0.8,2.5)
\end{pspicture}
\begin{pspicture}(0.0,0)(1.2,2.5)
 \psdots[dotsize=6pt](0.5,0.8)  \psdots[dotsize=6pt](0.5,1.8)
 \psline[arrows=->](0.5,0.8)(0.5,0.1) \rput(1.0,1.4){$\pm$}
 \psline[arrows=->](0.5,1.8)(0.5,0.9) \psline[arrows=->](0.2,1.5)(0.45,0.9)
 \psline[arrows=->](0.2,2.5)(0.45,1.9) \psline[arrows=->](0.5,1.8)(0.8,2.5)
\end{pspicture}
\begin{pspicture}(0.0,0)(1.4,2.5)
 \psdots[dotsize=6pt](0.5,0.8)  \psdots[dotsize=6pt](0.5,1.8)
 \psline[arrows=->](0.5,0.8)(0.5,0.1)   \rput(1.2,1.4){$\pm$}
 \psline[arrows=<-](0.5,1.7)(0.5,0.8)   \psline[arrows=<-](0.5,1.9)(0.5,2.6)
 \psline[arrows=->](0.15,1.5)(0.45,0.9) \psline[arrows=<-](0.8,1.5)(0.5,0.8)
\end{pspicture}
\begin{pspicture}(0.1,0)(1.2,2.5)
 \psdots[dotsize=6pt](0.5,0.8)  \psdots[dotsize=6pt](0.5,1.8)
 \psline[arrows=->](0.5,0.8)(0.5,0.1)  \rput(1.0,1.4){$\pm$}
 \psline[arrows=->](0.5,0.8)(0.5,1.7)  \psline[arrows=->](0.5,2.6)(0.5,1.9)
 \psline[arrows=->](0.2,2.5)(0.45,1.9) \psline[arrows=->](0.5,1.8)(0.8,2.5)
\end{pspicture}
\begin{pspicture}(0.1,0)(1.2,2.5)
 \psdots[dotsize=6pt](0.5,0.8)  \psdots[dotsize=6pt](0.5,1.8)
 \psline[arrows=->](0.5,0.8)(0.5,0.1)  \rput(1.0,1.4){$\pm$}
 \psline[arrows=<-](0.5,0.9)(0.5,1.7)  \psline[arrows=->](0.5,2.6)(0.5,1.9)
 \psline[arrows=->](0.2,2.5)(0.45,1.9) \psline[arrows=->](0.5,1.8)(0.8,2.5)
\end{pspicture}
\begin{pspicture}(0.1,0)(1.4,2.5)
 \psdots[dotsize=6pt](0.5,0.8)  \psdots[dotsize=6pt](0.5,1.8)
 \psline[arrows=->](0.5,0.8)(0.5,0.1) \rput(1.2,1.4){$\pm$}
 \psline[arrows=->](0.5,1.8)(0.5,0.9)  \psline[arrows=<-](0.8,1.5)(0.5,0.8)
 \psline[arrows=->](0.2,2.5)(0.45,1.9) \psline[arrows=->](0.8,2.5)(0.55,1.9)
\end{pspicture}
\begin{pspicture}(0.2,0)(1.9,2.5)
 \psdots[dotsize=6pt](0.5,0.8)  \psdots[dotsize=6pt](0.5,1.8)
 \psline[arrows=->](0.5,0.8)(0.5,0.1)   \rput(1.6,1.4){$=0$}
 \psline[arrows=->](0.5,1.7)(0.5,0.9)   \psline[arrows=<-](0.5,1.9)(0.5,2.6)
 \psline[arrows=->](0.85,1.5)(0.55,0.9) \psline[arrows=<-](1.1,1.3)(0.6,0.85)
\end{pspicture}
\] \caption{The boundary condition for $\Gamma=\Gamma(2,0)$}
 \label{boundary-condition}
\end{figure}
\end{itemize}
\end{defn}

\begin{remk}\label{weak-V-infty} A few comments are in order.
\begin{enumerate}

\item It is important to recall that the grading of $A^\ast$ is
negative of that of $A$.

\item The boundary condition interlocks operations of different
types.

\item In the above definition, we did not include $v_0 \in A$ only
to avid further consideration arising from the existence of a
homotopy unit for weak A$_\infty$-algebras. Nonetheless, one could
consider a discussion which allows for a $v_0 \in A$.
\end{enumerate}
\end{remk}

The concept of a V$_k$-algebra is homotopy invariant. That is, the
V$_k$-algebra structure on a complex can be transported onto a
quasi-isomorphic complex. A V$_1$-algebras is precisely an
A$_\infty$-algebra, and a V$_2$-algebra is exactly an
A$_\infty$-algebra with an invariant and symmetric homotopy
co-inner product. V$_3$-algebras are less familiar as they deal
with type $3$ vertices. An example of a vertex of higher type
appeared in H. Kajiura's thesis \cite[section 4.2]{K} in
connection with open strings.

\section{$\mathcal{DG}_\infty^\bullet$ and string
operations}\label{action section}

We now define a complex generated by graphs with vertices of
arbitrary types; see section \ref{vert O, I, II section}. We will
show how this complex acts on the Connes' cyclic Hochschild
complex. Intuitively, this generalizes the complex of cyclic
Sullivan chord diagrams by adding magnifying glass information to
the chords; see Figure \ref{strict}.

\subsection{Cyclic Hochschild complex}\label{CC section}

Let us recall the definition of the cyclic Hochschild complex. Let
$A$ be a graded vector space and define the cyclic Hochschild
complex $CC^\bullet(A)$ by,
\begin{multline*}
CC^\bullet(A)=\prod_{j=0}^\infty \Big\{f:A^{\otimes j+1}\to {\bf
k} \Big|
\\ f(a_1,\cdots,a_{j+1})=(-1)^{\epsilon} \cdot f(a_{j+1},a_1,\cdots, a_j)
\Big\},
\end{multline*}
where $\epsilon=|a_{j+1}|\cdot(|a_1|+\cdots+ |a_j|)$. Here,
$|\cdot|$ denotes the shifted degree of an element defined as
$|a|=deg(a)+1 $. A type zero vertex with six edges is shown in
Figure \ref{cyclic hochschild}.  We will later see how an element $f
\in (A^*)^{\otimes (j+1)} \subset CC^\bullet (A)$ with $j+1 \geq 6$
incoming edges can be placed on top of this vertex which plays the
role of an input. If the number of tensor factors, $j+1$, exceeds
$6$, then six of them are placed exactly on the top of the edges and
the remanning ones are left as what will be referred to as hairs
sticking out of the vertex.  This is done in all possible
combinatorial ways.
\begin{figure}
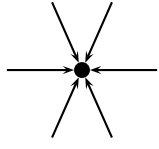

\[
    \psset{arrows=->}
    \pstree[treemode=L, levelsep=1cm, treesep=0.8cm]{\Tp}{
    \pstree[levelsep=0]{\Tc*{3pt}}{
    \pstree[levelsep=0.1cm,treemode=U,arrows=<-]{\Tn}{\Tp\Tp}
    \pstree[levelsep=1cm,treemode=L,arrows=<-]{\Tn}{\Tp}
    \pstree[levelsep=0.1cm,treemode=D,arrows=<-]{\Tn}{\Tp\Tp}
    }}
\] \caption{Cyclic Hochschild complex $CC^\bullet(A)$} \label{cyclic
hochschild}
\end{figure}

As it turns out, one does not need to specify a starting edge.
That is, if $f$ is a cyclic element, $f=\sum_{i=1}^{j+1}\pm a^*_i
\otimes \cdots \otimes a^*_{j+1}\otimes a^*_1\otimes \cdots
\otimes a^*_{i-1} \in (A^*)^{\otimes j+1}$, then we can choose any
of the $j+1$ incoming edges of Figure \ref{cyclic hochschild} as
the first edge, and insert the $a^*_1,\cdots, a^*_{j+1}$
cyclically. This procedure, as we shall see, is independent of the
choice of this first edge, due to cyclicity of the expressions and
an overall sign obtained from an orientation.

\begin{figure}
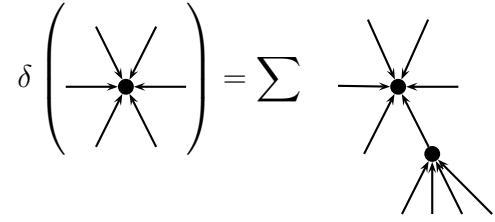

\[ \delta\left(
    \psset{arrows=->}
    \pstree[treemode=L, levelsep=0.8cm, treesep=0.8cm]{\Tp}{
    \pstree[levelsep=0]{\Tc*{3pt}}{
    \pstree[levelsep=0.08cm,treemode=D,arrows=<-]{\Tn}{\Tp\Tp}
    \pstree[levelsep=0.8cm,treemode=L,arrows=<-]{\Tn}{\Tp}
    \pstree[levelsep=0.08cm,treemode=U,arrows=<-]{\Tn}{\Tp\Tp}
    }}\right)=\sum \quad\psset{arrows=->}
    \pstree[treemode=L, levelsep=0.8cm, treesep=0.8cm]{\Tp}{
    \pstree[levelsep=0]{\Tc*{3pt}}{
    \pstree[levelsep=0.08cm,treemode=U,arrows=<-]{\Tn}{\Tp\Tp}
    \pstree[levelsep=0.8cm,treemode=L,arrows=<-]{\Tn}{\Tp}
    \pstree[levelsep=0.08cm,treemode=D,arrows=<-]{\Tn}{\Tp
    \pstree[levelsep=0.8cm,treesep=0.4cm]{\Tc*{3pt}}{\Tn\Tp\Tp\Tp\Tp}}
    }} \] \caption{Differential on $CC^\bullet(A)$}
\label{differential on CC(A)}
\end{figure}

Now, for an A$_\infty$-algebra $A$, one can define a differential
$\delta:CC^\bullet(A)\to CC^\bullet(A)$ as pictured in Figure
\ref{differential on CC(A)}.
To be more precise, the differential $\delta$ applies the
A$_\infty$-algebra maps $\mu_m$ (for $m\geq 1$) in all possible
ways to an element $f\in CC^\bullet(A)$, in the following fashion,
\begin{eqnarray*}
\delta f (a_1,\cdots,a_{j+1})&=& \quad\sum_{l,i} \pm f\big(
\mu_l(a_{i}, \cdots, a_{j+1},a_1,\cdots),\cdots, a_{i-1}\big)\\ &&+
\sum_{l,i} \pm f \big(a_1, \cdots, \mu_l(a_i,\cdots, a_{i+l-1}),
\cdots , a_{j+1} \big)\\&&+ \sum_{l,i} \pm f \big( a_i,\cdots, \mu_l
(\cdots, a_{j+1} ,a_1, \cdots, a_{i-1})\big).\label{diff-formula}
\end{eqnarray*}
The signs are determined by comparing the linear order of the
symbols in an expression to the orientation $\mu_l\wedge f \wedge
a_1\wedge \cdots\wedge a_{j+1}$, where $\mu_l$ is of degree $1$,
and $a_i$ has the shifted degree $|a_i|=deg(a_i)+1$.  Note that if
$f\in (A^*)^{\otimes j+1}$ is cyclically invariant, then
$\delta(f)$ is also cyclically invariant. Furthermore, the picture
for $\delta^2(f)$ corresponds to applying two A$_\infty$-algebra
maps to $f$, say $\mu_l$ and $\mu_m$. Note that each term for
which $\mu_l$ and $\mu_m$ sit on different edges of $f$ appears
twice, obtained from first applying $\mu_l$ and then $\mu_m$, and
vice versa. Since $\mu_m$ moves over $\mu_l$ in one case, but not
the other, they appear with opposite signs and cancel. As for the
terms where $\mu_l$ sits on top of $\mu_m$, they add up to zero
due to the compatibility requirements of the $\mu_i$'s in the
definition of an A$_\infty$-algebra; see section
\ref{A-infty-section}. This shows that $\delta^2=0$. The homology
of this complex, known as the cyclic Hochschild cohomology, is
denoted by $HC^\bullet(A)$.

\subsection{PROP $\mathcal{DG}_k^\bullet$ of directed graphs of type $k$}\label{graph complex section}

\begin{figure}
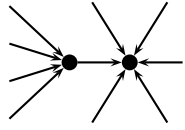

\[
    \psset{arrows=->}
    \pstree[treemode=L, levelsep=0.7cm]{\Tp}{
    \pstree[levelsep=0, arrows=<-, treesep=0.5cm]{\Tc*{3pt}}{
    \pstree[levelsep=0.3cm,treemode=U]{\Tn}{\Tp\Tn\Tp}
    \pstree[levelsep=0.8cm,treemode=L]{\Tn}{\pstree{\Tc*{3pt}}{\Tp\Tp\Tp\Tp}}
    \pstree[levelsep=0.3cm,treemode=D]{\Tn}{\Tp\Tn\Tp}
    }}
 \] \caption{Plugging type $1$ to type $0$}
\label{O in I}
\end{figure}

One basic idea is that an incoming edge can connect to an outgoing
edge, as indicated in Figure \ref{O in I}.
We are interested in graphs that are obtained by connecting
vertices of mixed types.
\begin{defn} $\mathcal{DG}_\infty^\bullet$ is the vector space generated by
all pairs $(\gamma, \mathcal O)$, where $\gamma$ is a directed graph
with a cyclic ordering of the edges at each vertex of type
$0,1,2,\cdots$, with the proviso that type $1$ vertices have at
least two incoming edges, and $\mathcal O$ is an orientation on the
graded vector space generated by the edges and vertices of type
$n\neq 0$ of $\gamma$, where the degree of an edge is $1$, and the
degree of a vertex of type $n$ is $n(2-d)+(d-4)$.

Furthermore, we define $\mathcal{DG}_k^\bullet$, for $k=1,2,3
\cdots$, to be the subspace of $\mathcal{DG}_\infty^\bullet$ given
by graphs with vertices of type $n$ where $n \leq k$.
\end{defn}
\begin{figure}
\[
\begin{pspicture}(0,0)(6,4)
 \psdots[dotsize=6pt](1.4,3.2) \psdots[dotsize=6pt](4.6,3.2)
 \psdots[dotsize=6pt](0.8,0.8) \psdots[dotsize=6pt](1.8,2)
 \psdots[dotsize=6pt](2.6,2) \psdots[dotsize=6pt](3.4,2)
 \psdots[dotsize=6pt](4,0.4) \psdots[dotsize=6pt](5.2,1)
 \psline[arrows=->](2.2,2)(1.9,2) \psline[arrows=->](2.2,2)(3.3,2)
 \psline[arrows=->](0.8,0.8)(1.8,1.9) \psline[arrows=->](0.8,0.8)(1.35,3.1)
 \psline[arrows=->](4,0.4)(1.88,1.92) \psline[arrows=->](4,0.4)(4.55,3.1)
 \psline[arrows=->](3.4,2)(1.45,3.1) \psline[arrows=->](3.4,2)(3.95,0.5)
 \psline[arrows=->](5.2,1)(3.5,2) \psline[arrows=->](5.2,1)(4.65,3.1)
 \psline[arrows=->](1.8,2)(4.5,3.2)
\end{pspicture}
\]
\caption{A combinatorial directed graph} \label{a-graph}
\end{figure}
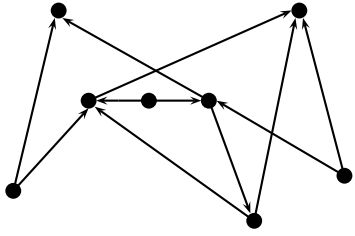

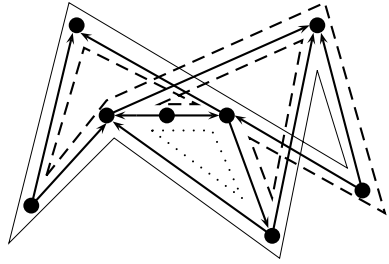
\begin{figure}
\[
\begin{pspicture}(0,0)(6,4)
 \psdots[dotsize=6pt](1.4,3.2) \psdots[dotsize=6pt](4.6,3.2)
 \psdots[dotsize=6pt](0.8,0.8) \psdots[dotsize=6pt](1.8,2)
 \psdots[dotsize=6pt](2.6,2) \psdots[dotsize=6pt](3.4,2)
 \psdots[dotsize=6pt](4,0.4) \psdots[dotsize=6pt](5.2,1)
 \psline[arrows=->](2.2,2)(1.9,2) \psline[arrows=->](2.2,2)(3.3,2)
 \psline[arrows=->](0.8,0.8)(1.8,1.9) \psline[arrows=->](0.8,0.8)(1.35,3.1)
 \psline[arrows=->](4,0.4)(1.88,1.92) \psline[arrows=->](4,0.4)(4.55,3.1)
 \psline[arrows=->](3.4,2)(1.45,3.1) \psline[arrows=->](3.4,2)(3.95,0.5)
 \psline[arrows=->](5.2,1)(3.5,2) \psline[arrows=->](5.2,1)(4.65,3.1)
 \psline[arrows=->](1.8,2)(4.5,3.2)
 \psline[linestyle=dotted](2.4,1.8)(3.2,1.8)(3.6,0.9)(2.4,1.8)
 \psline[linestyle=solid, linewidth=0.05pt](4.6,2.6)(5,1.3)(1.3,3.5)(0.5,0.3)(1.9,1.7)(4.1,0.1)(4.6,2.6)
 \psline[linestyle=dashed](1,1.2)(1.8,2.2)(4.7,3.5)(5.5,0.7)(3.7,1.7)(4,0.9)
 (4.4,3)(2.5,2.15)(3,2.15)(1.5,2.9)(1,1.2)
\end{pspicture}
\]
\caption{2 inputs and 3 outputs} \label{outputs}
\end{figure}

Figure \ref{a-graph} shows an example of a graph in
$\mathcal{DG}_\infty^\bullet$. For each such graphs, we have so
called inputs and outputs. By definition, the inputs are the
vertices of type $0$. The outputs of the graph are the boundaries of
the thickened (fat) graph. See for example Figure \ref{outputs},
where the three outputs are indicated by dotted, dashed, and thin
solid curves. Both inputs and outputs are enumerated as part of the
structure. Note that there may be no input, but there has to be at
least one output.

\begin{figure}
\[
\begin{pspicture}(0,0)(7,5)
 \psdots[dotsize=6pt](0.4,4.6) \psdots[dotsize=6pt](3,4) \psdots[dotsize=6pt](6,4.4)
 \psdots[dotsize=6pt](6.6,2.4) \psdots[dotsize=6pt](6.6,1.4) \psdots[dotsize=6pt](6.6,0.4)
 \psdots[dotsize=6pt](0.6,1) \psdots[dotsize=6pt](1.8,1.4) \psdots[dotsize=6pt](1,3)
 \psdots[dotsize=6pt](3.2,2.4) \psdots[dotsize=6pt](3,1.6) \psdots[dotsize=6pt](4.2,1.2)
 \psdots[dotsize=6pt](5.6,3.4) \psdots[dotsize=6pt](4.8,2.8) \psdots[dotsize=6pt](5.2,1)
 \psdots[dotsize=6pt](5.6,0)
 \psline[arrows=->](1.8,1.4)(0.7,1) \psline[arrows=->](1,3)(1.73,1.47)
 \psline[arrows=->](6.6,1.4)(6.6,2.3) \psline[arrows=->](6.6,1.4)(6.6,0.5)
 \psline[arrows=->](3.2,2.4)(3,1.7) \psline[arrows=->](3,1.6)(1.9,1.4)
 \psline[arrows=->](4.2,1.2)(3.1,1.6) \psline[arrows=->](4.8,2.8)(4.2,1.3)
 \psline[arrows=->](5.2,1)(4.3,1.2) \psline[arrows=->](5.2,1)(5.6,0.1)
 \psline[arrows=->](5.6,3.4)(6,4.3) \psline[arrows=->](5.6,3.4)(4.88,2.88)
 \psline[arrows=->](3,4)(4.72,2.88)
 \pscurve[arrows=->](3,4)(3,3)(3.2,2.5) \pscurve[arrows=->](1,3)(3.5,2.9)(3.3,2.45)
\end{pspicture}
\]
\caption{A graph in $\mathcal{DG}_\infty^0$} \label{lowest graph}
\end{figure}
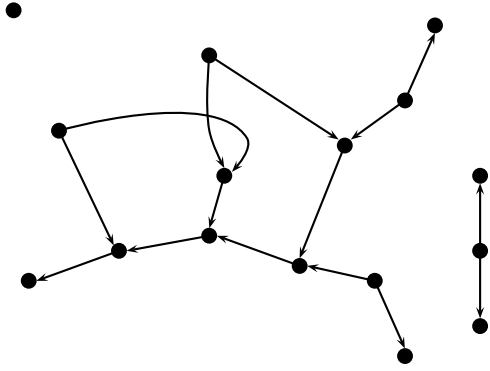

The grading in $\mathcal{DG}_\infty^\bullet$ is given in the
following way. Look at a graph in $\mathcal{DG}_\infty^\bullet$.
Every nontrivial vertex of type $0$ with $m$ incoming edges
contributes $m-1$ to the degree, whereas vertices of type $n\geq 1$
with $m$ incoming edges add $2n+ m-2$ to the degree. Trivial
vertices have degree zero. This introduces a grading
$\mathcal{DG}_\infty^\bullet= \bigoplus_{ n\in \mathbb N \cup \{0\}}
\mathcal{DG}_\infty^n$. Note that none of the graphs have negative
degree in $\mathcal{DG}_\infty^\bullet$, since we restricted type
$1$ vertices to have at least two incoming edges. In fact, the
lowest degree $\mathcal{DG}_\infty^0$ consists of type $0$ vertices
with zero or one incoming edge, type $1$ vertices with two incoming
edges, and type $2$ vertices with no incoming edges; see Figure
\ref{lowest graph}.
Thus, the lowest operations are given by the product and the
co-inner product, whereas all higher vertices, including all
vertices of type $\geq 3$, are interpreted as homotopies.

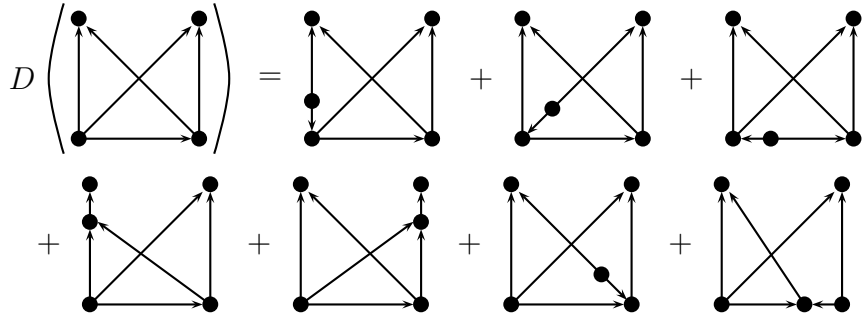
\begin{figure}
\[
\begin{pspicture}(-0.8,0)(2.3,2)
 \psdots[dotsize=6pt](0.2,0.2)  \psdots[dotsize=6pt](1.8,0.2)
 \psdots[dotsize=6pt](0.2,1.8)  \psdots [dotsize=6pt](1.8,1.8)
 \psline[arrows=->](0.2,0.2)(0.2,1.7) \psline[arrows=->](1.8,0.2)(1.8,1.7)
 \psline[arrows=->](0.2,0.2)(1.7,1.7) \psline[arrows=->](1.8,0.2)(0.3,1.7)
 \psline[arrows=->](0.2,0.2)(1.7,0.2) \rput(-0.55,1){$D$}
 \pscurve(0,0)(-0.2,1)(0,2) \pscurve(2,0)(2.2,1)(2,2)
\end{pspicture}
\begin{pspicture}(-0.8,0)(2,2)
 \psdots[dotsize=6pt](0.2,0.2)  \psdots[dotsize=6pt](1.8,0.2)
 \psdots[dotsize=6pt](0.2,1.8)  \psdots[dotsize=6pt](1.8,1.8)
 \psline[arrows=->](0.2,0.7)(0.2,1.7)
 \psline[arrows=->](0.2,0.7)(0.2,0.3) \psdots[dotsize=6pt](0.2,0.7)
 \psline[arrows=->](1.8,0.2)(1.8,1.7)
 \psline[arrows=->](0.2,0.2)(1.7,1.7) \psline[arrows=->](1.8,0.2)(0.3,1.7)
 \psline[arrows=->](0.2,0.2)(1.7,0.2) \rput(-0.35,1){$=$}
\end{pspicture}
\begin{pspicture}(-0.8,0)(2,2)
 \psdots[dotsize=6pt](0.2,0.2)  \psdots[dotsize=6pt](1.8,0.2)
 \psdots[dotsize=6pt](0.2,1.8)  \psdots[dotsize=6pt](1.8,1.8)
 \psdots[dotsize=6pt](0.6,0.6) \psline[arrows=->](0.6,0.6)(0.26,0.26)
 \psline[arrows=->](0.2,0.2)(0.2,1.7) \psline[arrows=->](1.8,0.2)(1.8,1.7)
 \psline[arrows=->](0.6,0.6)(1.7,1.7) \psline[arrows=->](1.8,0.2)(0.3,1.7)
 \psline[arrows=->](0.2,0.2)(1.7,0.2) \rput(-0.35,1){$+$}
\end{pspicture}
\begin{pspicture}(-0.8,0)(2,2)
 \psdots[dotsize=6pt](0.2,0.2)  \psdots[dotsize=6pt](1.8,0.2)
 \psdots[dotsize=6pt](0.2,1.8)  \psdots[dotsize=6pt](1.8,1.8)
 \psline[arrows=->](0.2,0.2)(0.2,1.7) \psline[arrows=->](1.8,0.2)(1.8,1.7)
 \psline[arrows=->](0.2,0.2)(1.7,1.7) \psline[arrows=->](1.8,0.2)(0.3,1.7)
 \psline[arrows=->](0.7,0.2)(1.7,0.2) \rput(-0.35,1){$+$}
 \psdots[dotsize=6pt](0.7,0.2) \psline[arrows=->](0.7,0.2)(0.3,0.2)
\end{pspicture}
\]
\[
\begin{pspicture}(-0.8,0)(2,2)
 \psdots[dotsize=6pt](0.2,0.2)  \psdots[dotsize=6pt](1.8,0.2)
 \psdots[dotsize=6pt](0.2,1.8)  \psdots[dotsize=6pt](1.8,1.8)
 \psline[arrows=->](0.2,0.2)(0.2,1.2) \psline[arrows=->](1.8,0.2)(1.8,1.7)
 \psline[arrows=->](0.2,0.2)(1.7,1.7) \psline[arrows=->](1.8,0.2)(0.3,1.3)
 \psline[arrows=->](0.2,0.2)(1.7,0.2) \rput(-0.35,1){$+$}
 \psdots[dotsize=6pt](0.2,1.3) \psline[arrows=->](0.2,1.3)(0.2,1.7)
\end{pspicture}
\begin{pspicture}(-0.8,0)(2,2)
 \psdots[dotsize=6pt](0.2,0.2)  \psdots[dotsize=6pt](1.8,0.2)
 \psdots[dotsize=6pt](0.2,1.8)  \psdots[dotsize=6pt](1.8,1.8)
 \psline[arrows=->](0.2,0.2)(0.2,1.7) \psline[arrows=->](1.8,0.2)(1.8,1.2)
 \psline[arrows=->](0.2,0.2)(1.7,1.3) \psline[arrows=->](1.8,0.2)(0.3,1.7)
 \psline[arrows=->](0.2,0.2)(1.7,0.2) \rput(-0.35,1){$+$}
 \psdots[dotsize=6pt](1.8,1.3) \psline[arrows=->](1.8,1.3)(1.8,1.7)
\end{pspicture}
\begin{pspicture}(-0.8,0)(2,2)
 \psdots[dotsize=6pt](0.2,0.2)  \psdots[dotsize=6pt](1.8,0.2)
 \psdots[dotsize=6pt](0.2,1.8)  \psdots[dotsize=6pt](1.8,1.8)
 \psline[arrows=->](0.2,0.2)(0.2,1.7) \psline[arrows=->](1.8,0.2)(1.8,1.7)
 \psline[arrows=->](0.2,0.2)(1.7,1.7) \psline[arrows=->](1.4,0.6)(0.3,1.7)
 \psline[arrows=->](0.2,0.2)(1.7,0.2) \rput(-0.35,1){$+$}
 \psdots[dotsize=6pt](1.4,0.6) \psline[arrows=->](1.4,0.6)(1.74,0.26)
\end{pspicture}
\begin{pspicture}(-0.8,0)(2,2)
 \psdots[dotsize=6pt](0.2,0.2)  \psdots[dotsize=6pt](1.8,0.2)
 \psdots[dotsize=6pt](0.2,1.8)  \psdots[dotsize=6pt](1.8,1.8)
 \psline[arrows=->](0.2,0.2)(0.2,1.7) \psline[arrows=->](1.8,0.2)(1.8,1.7)
 \psline[arrows=->](0.2,0.2)(1.7,1.7) \psline[arrows=->](1.3,0.2)(0.3,1.7)
 \psline[arrows=->](0.2,0.2)(1.2,0.2) \rput(-0.35,1){$+$}
 \psdots[dotsize=6pt](1.3,0.2) \psline[arrows=<-](1.4,0.2)(1.8,0.2)
\end{pspicture}
\] \caption{Boundary of an element in $\mathcal{DG}_\infty^\bullet$} \label{boundary-D}
\end{figure}

$\mathcal{DG}_\infty^\bullet$ is in fact a complex. The differential
$D: \mathcal{DG}_\infty^\bullet \to \mathcal{DG}_\infty^{\bullet-
1}$ is defined by expansion of an edge. To be more precise, $D$
applied to a graph is defined to be the sum of all permitted
directed graphs that preserve the number of type $0$ vertices, such
that collapsing an edge produces the original graph. An example is
shown in Figure \ref{boundary-D}, where the first three graphs of
the boundary expand the lower left type $3$ vertex, the next two
graphs expand the upper two vertices of type $0$, and the last two
graphs expand the lower right vertex of type $2$.
The orientations of the resulting graphs, which we suppressed in
Figure \ref{boundary-D}, are obtained by bringing the relevant
vertex in the expression determining the orientation to the front,
and replacing it by the added edge followed by the two vertices in
the order given by the direction of the expanded edge. Note that
$D^2=0$, since $D^2$ expands two edges, which can be done in two
ways, depending on the ordering of the expanded edges.

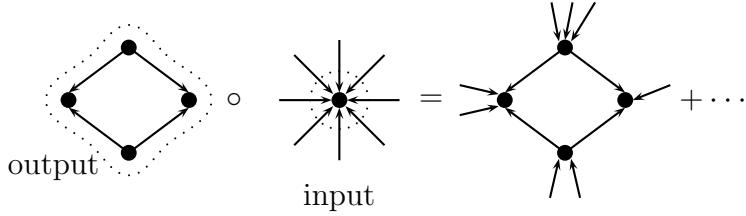
\begin{figure}
\[
\begin{pspicture}(0.4,0)(3.8,3)
 \psdots[dotsize=6pt](2,2.2) \psdots[dotsize=6pt](2,0.8)
 \psdots[dotsize=6pt](1.2,1.5) \psdots[dotsize=6pt](2.8,1.5)
 \psline[arrows=->](2,2.2)(1.2,1.6) \psline[arrows=->](2,2.2)(2.8,1.6)
 \psline[arrows=->](2,0.8)(1.2,1.4) \psline[arrows=->](2,0.8)(2.8,1.4)
 \psccurve[linestyle=dotted](2,2.5)(2.5,2.1)(3.1,1.5)(2.5,0.9)(2,0.5)(1.5,0.9)(0.9,1.5)(1.5,2.1)
 \rput(1,0.6){output} \rput(3.4,1.5){$\circ$}
\end{pspicture}
\begin{pspicture}(-2,0)(4.4,3)
 \psdots[dotsize=6pt](2,2.2) \psdots[dotsize=6pt](2,0.8)
 \psdots[dotsize=6pt](-1,1.5)
 \psdots[dotsize=6pt](1.2,1.5) \psdots[dotsize=6pt](2.8,1.5)
 \psline[arrows=->](2,2.2)(1.2,1.6) \psline[arrows=->](2,2.2)(2.8,1.6)
 \psline[arrows=->](2,0.8)(1.2,1.4) \psline[arrows=->](2,0.8)(2.8,1.4)
 \psline[arrows=->](1.82,2.8)(1.92,2.28) \psline[arrows=->](2.1,2.8)(2,2.3)
 \psline[arrows=->](2.4,2.8)(2.08,2.28)
 \psline[arrows=->](1.8,0.2)(1.92,0.72) \psline[arrows=->](2.2,0.2)(2.08,0.72)
 \rput(4,1.5){$+\cdots$} \rput(0.22,1.5){$=$} \rput(-1,0.2){input}
 \psline[arrows=->](-1,0.7)(-1,1.4) \psline[arrows=->](-1,2.3)(-1,1.6)
 \psline[arrows=->](-1.8,1.5)(-1.1,1.5) \psline[arrows=->](-0.2,1.5)(-0.9,1.5)
 \psline[arrows=->](-1.6,2.1)(-1.08,1.58) \psline[arrows=->](-0.4,2.1)(-0.92,1.58)
 \psline[arrows=->](-1.6,0.9)(-1.08,1.42) \psline[arrows=->](-0.4,0.9)(-0.92,1.42)
 \pscircle[linestyle=dotted](-1,1.5){0.4}
 \psline[arrows=->](3.4,1.7)(2.9,1.5)
 \psline[arrows=->](0.6,1.7)(1.12,1.58) \psline[arrows=->](0.6,1.3)(1.12,1.42)
\end{pspicture}
\]\caption{The composition} \label{composition}
\end{figure}

We next explain how to compose two graphs. Suppose $\gamma_1$ is a
graph with $i$ inputs and $j$ outputs, and let $\gamma_2$ be a
graph with $j$ inputs and $k$ outputs. The composition will then
be a graph $\gamma_1\circ \gamma_2$ with $i$ inputs and $k$
outputs. Identify the outputs of $\gamma_1$ with the inputs of
$\gamma_2$. Then, starting with the graph $\gamma_1$, we perform
the operation from Figure \ref{composition} at each output of
$\gamma_1$.
More precisely, we take the edges at an input of $\gamma_2$ and
attach them to the corresponding output of $\gamma_1$ in all
possible ways. The orientation of the composed graph
$\gamma_1\circ \gamma_2$ is given by attaching the orientation of
$\gamma_2$ to that of $\gamma_1$ on the right. Note that if we
thicken the graphs to a surface and punch holes at locations of
the type zero vertices, the composition corresponds exactly to the
gluing of the surfaces along $j$ boundary components. This picture
also shows that the composed graph has in fact exactly $k$
outputs.

The above definitions give $\mathcal{DG}_\infty^\bullet$ and each
$\mathcal{DG}_k^\bullet$ the structure of a PROP. Note that there
are obvious inclusion maps $\mathcal{DG}_k^\bullet \hookrightarrow
\mathcal{DG}_{k+1}^\bullet$ and $\mathcal{DG}_k^\bullet
\hookrightarrow \mathcal{DG}_\infty^\bullet$, for every $k=1,2,
\cdots$.
\begin{prop}
There is a sequence of inclusions of PROPs, $$
\mathcal{DG}_1^\bullet \hookrightarrow \mathcal{DG}_2^\bullet
\hookrightarrow \mathcal{DG}_3^\bullet \hookrightarrow \cdots
\hookrightarrow \mathcal{DG}_\infty^\bullet. $$
\end{prop}

\subsection{Algebraic structure of the cyclic Hochschild complex}
\label{action on CC section}

We are now prepared to describe how, for an V$_\infty$-algebra
$A$, the PROP $\mathcal{DG}_\infty^\bullet$ acts on
$CC^\bullet(A)$.
\begin{thrm}\label{V_k-action}
Let $k=1,2,\cdots,\infty$, and let $A$ be a V$_k$-algebra. Then,
there is a map, $$ \mathcal{DG}_k^\bullet \to \prod_{r,s}Hom\Big(
CC^\bullet(A)^{ \otimes r}, CC^\bullet( A)^{ \otimes s}\Big),$$
which respects the grading, differential, and composition.
\end{thrm}

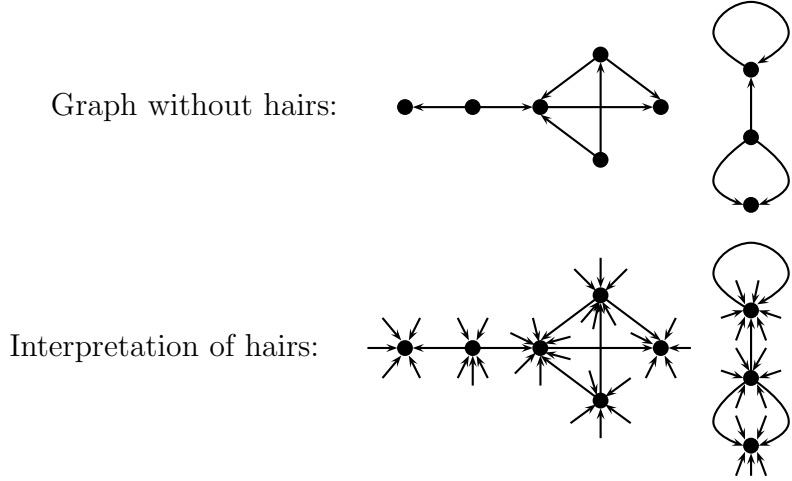
\begin{figure}
\[
\begin{pspicture}(-6.2,0)(4.6,3)
 \rput(-3.4,1.5){Graph without hairs:}
 \psdots[dotsize=6pt](4,0.2) \psdots[dotsize=6pt](4,1.1) \psdots[dotsize=6pt](4,2)
 \psdots[dotsize=6pt](2,2.2) \psdots[dotsize=6pt](2,0.8)
 \psdots[dotsize=6pt](1.2,1.5) \psdots[dotsize=6pt](2.8,1.5)
 \psdots[dotsize=6pt](0.3,1.5) \psdots[dotsize=6pt](-0.6,1.5)
 \psline[arrows=->](2,2.2)(1.2,1.6) \psline[arrows=->](2,2.2)(2.8,1.6)
 \psline[arrows=->](2,0.8)(1.2,1.4) \psline[arrows=->](1.2,1.5)(2.7,1.5)
 \psline[arrows=->](0.3,1.5)(-0.5,1.5) \psline[arrows=->](0.3,1.5)(1.1,1.5)
 \psline[arrows=->](2,0.8)(2,2.1) \psline[arrows=->](4,1.1)(4,1.9)
 \pscurve[arrows=->](4,1.1)(3.5,0.5)(3.9,0.2)
 \pscurve[arrows=->](4,1.1)(4.5,0.5)(4.1,0.2)
 \pscurve[arrows=->](4,2)(3.5,2.5)(4,2.9)(4.5,2.5)(4.08,2.08)
\end{pspicture}
\]
\[
\begin{pspicture}(-6.2,0)(4.6,3)
 \rput(-3.8,1.5){Interpretation of hairs:}
 \psdots[dotsize=6pt](2,0.8)
 \psline[arrows=->](1.85,1.2)(1.95,0.87) \psline[arrows=->](2,0.3)(2,0.7)
 \psline[arrows=->](2.4,0.5)(2.06,0.75) \psline[arrows=->](2.4,1.1)(2.06,0.85)
 \psline[arrows=->](1.6,0.5)(1.94,0.75)
 \psdots[dotsize=6pt](1.2,1.5)
 \psline[arrows=->](1.1,1.9)(1.18,1.58) \psline[arrows=->](1.2,1)(1.2,1.4)
 \psline[arrows=->](0.86,1.16)(1.14,1.44) \psline[arrows=->](1.6,1.65)(1.26,1.56)
 \psline[arrows=->](1.6,1.35)(1.26,1.44) \psline[arrows=->](0.8,1.7)(1.13,1.54)
 \psdots[dotsize=6pt](4,0.2)
 \psline[arrows=->](3.8,0.6)(3.92,0.28) \psline[arrows=->](4.2,0.6)(4.08,0.28)
 \psline[arrows=->](4.2,-0.2)(4.08,0.12) \psline[arrows=->](4,-0.2)(4,0.1)
 \psline[arrows=->](3.8,-0.2)(3.92,0.12)
 \psdots[dotsize=6pt](4,1.1)
 \psline[arrows=->](3.8,0.7)(3.92,1.02) \psline[arrows=->](4.2,0.7)(4.08,1.02)
 \psline[arrows=->](4.2,1.5)(4.02,1.19) \psline[arrows=->](4.4,1.2)(4.1,1.1)
 \psline[arrows=->](3.8,1.5)(3.98,1.19) \psline[arrows=->](3.6,1.2)(3.9,1.1)
 \psdots[dotsize=6pt](4,2)
 \psline[arrows=->](3.8,2.4)(3.92,2.08) \psline[arrows=->](4.1,2.4)(4.03,2.08)
 \psline[arrows=->](4.2,1.6)(4.08,1.92) \psline[arrows=->](3.8,1.6)(3.92,1.92)
 \psline[arrows=->](4.4,1.9)(4.1,2) \psline[arrows=->](3.6,1.9)(3.9,2)
 \psdots[dotsize=6pt](-0.6,1.5)
 \psline[arrows=->](-0.9,1.9)(-0.65,1.6) \psline[arrows=->](-0.4,1.9)(-0.55,1.6)
 \psline[arrows=->](-0.9,1.1)(-0.65,1.4) \psline[arrows=->](-0.4,1.1)(-0.55,1.4)
 \psline[arrows=->](-1.1,1.5)(-0.7,1.5)
 \psdots[dotsize=6pt](2.8,1.5)
 \psline[arrows=->](2.4,1.7)(2.7,1.55)
 \psline[arrows=->](2.6,1.1)(2.75,1.4) \psline[arrows=->](3,1.1)(2.85,1.4)
 \psline[arrows=->](3.2,1.5)(2.9,1.5)  \psline[arrows=->](3,1.9)(2.85,1.6)
 \psdots[dotsize=6pt](0.3,1.5)
 \psline[arrows=->](0.1,1.9)(0.25,1.6) \psline[arrows=->](0.5,1.9)(0.35,1.6)
 \psline[arrows=->](0.1,1.1)(0.25,1.4) \psline[arrows=->](0.5,1.1)(0.35,1.4)
 \psline[arrows=->](0.3,1)(0.3,1.4)
 \psdots[dotsize=6pt](2,2.2)
 \psline[arrows=->](2,2.7)(2,2.3) \psline[arrows=->](2.2,1.8)(2.02,2.12)
 \psline[arrows=->](1.9,1.75)(1.98,2.12) \psline[arrows=->](1.74,1.84)(1.98,2.12)
 \psline[arrows=->](2.35,2.55)(2.06,2.26) \psline[arrows=->](1.65,2.55)(1.94,2.26)
 \psline[arrows=->](2,2.2)(1.2,1.6) \psline[arrows=->](2,2.2)(2.8,1.6)
 \psline[arrows=->](2,0.8)(1.2,1.4) \psline[arrows=->](1.2,1.5)(2.7,1.5)
 \psline[arrows=->](0.3,1.5)(-0.5,1.5) \psline[arrows=->](0.3,1.5)(1.1,1.5)
 \psline[arrows=->](2,0.8)(2,2.1) \psline[arrows=->](4,1.1)(4,1.9)
 \pscurve[arrows=->](4,1.1)(3.5,0.5)(3.9,0.2)
 \pscurve[arrows=->](4,1.1)(4.5,0.5)(4.1,0.2)
 \pscurve[arrows=->](4,2)(3.5,2.5)(4,2.9)(4.5,2.5)(4.08,2.08)
\end{pspicture}
\]
\caption{Graphs interpreted with hairs} \label{typical graph}
\end{figure}
Here is a sketch of the proof for $\mathcal{DG}_\infty^\bullet$.
The proof for $\mathcal{DG}_k^\bullet$ is similar. For any graph
in $\mathcal{DG}_\infty^\bullet$ with $r$ inputs and $s$ outputs,
we will obtain a mapping $CC^\bullet(A)^{ \otimes r}\to
CC^\bullet( A)^{\otimes s}$. The elements $f_1,\cdots, f_r\in
CC^\bullet(A)$ are to be plugged into the type $0$ vertices
enumerated by $1,\cdots, r$. In detail, we consider the $j^{th}$
vertex of type $0$, together with all possible hairs, as indicated
in Figure \ref{typical graph}. If $f_j$ is given by a sum of terms
$a^*_1\otimes\cdots \otimes a^*_n$, then choose one of the edges
as your initial edge, and apply $a^*_1, \cdots, a^*_n$ cyclically
at this vertex. We shall see that the final outcome will be
independent of the choice of the starting point; see section
\ref{CC section}. Note that the total number of edges and hairs
must add up exactly to $n$, the number of the tensor factors. If
there are fewer than $n$ edges, then add an appropriate number of
hairs in all possible ways. If there are more than $n$ edges, then
$a^*_1,\cdots,a^*_n$ cannot be applied to the vertex and would
contribute a zero to the sum.

Next, we also apply the V$_\infty$-algebra structure to the
vertices of type $k\geq 1$. Just as with vertices of type $0$, we
always add the appropriate number of hairs and sum over all
possible ways. Since type $1$ vertices have a unique outgoing
edge, the application of elements of $A$ and $A^*$ is given
without ambiguity. Similarly, we assign the element
$v_{i_1,\cdots, i_k}$ of the V$_\infty$-algebra structure to a
vertices of type $k$. The symmetry condition from Definition
\ref{V_k} guarantees that this application is, up to a sign,
independent of choice of first outgoing edge; see Figure
\ref{homotopy co-inner product}. There is an overall sign
determined by the orientation $\mathcal O$ which fixes this
ambiguity in sign.

After applying these algebraic data to all vertices, using the
paring between $A^*$ and $A$, the dual elements on each of the
edges are contracted so that the only remaining terms are the
hairs. We can now read off the outputs along the outputs of the
graph. These outputs have no starting points and need to be
symmetrized cyclically to become elements of $CC^\bullet (A)$.

Note that the above procedure is only well defined element up to
sign. For one, the symmetry condition determines
$v_{i_1,\cdots,i_k}$ at a vertex only up to sign, and a second
problem arises when contracting $A^*$ with $A$, since we shifted
$A^*$, but not $A$. These ambiguities may be resolved, using the
orientation on the vertices and edges of the graph in
$\mathcal{DG}_\infty^\bullet$. In fact, using the orientation, as
well as the cyclic ordering of the edges at each vertex, we may
linearize the expression for all vertices and put the inverse of the
shift operators in those places where edges appear. The inverse of
the shift operators, one for each edge in the expression for the
orientation, will be used to shift back elements of $A^\ast$ before
evaluating them on elements of $A$. The remaining terms in the
linearization may be rearranged in accordance with the rules
dictated by the graph. The symmetry condition guarantees a well
defined overall sign for the outcome of the operation after the
rearrangement, since we only made repeated use of the Koszul rule.

In order to check that the differentials are preserved, recall that
the differential on $CC^\bullet(A)$ is given by applying one
additional A$_\infty$-algebra map (type $1$ vertex) as described in
Figure \ref{differential on CC(A)}. The differential in the space of
operations $Hom\big( CC^\bullet(A)^{ \otimes r}, CC^\bullet(A
)^{\otimes s}\big)$ applies this differential to all inputs and
outputs. This cancels all additional A$_\infty$-algebra maps on
hairs of type $0$ vertices, and adds a type $1$ vertex at the type
$0$ vertices in all possible ways. The remaining additional
A$_\infty$-algebra maps together with the boundary condition from
Definition \ref{V_k} shows that, in fact, one has to add one vertex
at any possible position of the graph. This coincides with the
differential in $\mathcal{DG}_\infty^\bullet$.

Finally, we deal with the compositions. Consider the usual
composition of $Hom\big( CC^\bullet(A)^{ \otimes r}, CC^\bullet(A
)^{\otimes s}\big)$ and $Hom\big( CC^\bullet(A)^{ \otimes s},
CC^\bullet(A )^{\otimes t}\big)$ with output in $Hom\big(
CC^\bullet (A)^{ \otimes r}, CC^\bullet(A )^{\otimes t}\big)$. In
terms of graphs in $\mathcal{DG}_\infty^\bullet$, this means that
the $s$ outputs of the first graph have to be taken as the inputs
for the second graph. We defined the $s$ outputs as the leftover
hairs after applying the edges of the graph. Thus, we need to
apply the edges at the inputs (type $0$ vertices) of the second
graph, to the outputs of the first graph in all possible ways
permitted by the combinatorics of the $s$ outputs. This is exactly
the definition of the composition as explained in Figure
\ref{composition}.

\subsection{An example: associative algebras with co-inner products}
\label{strict section}

Recall that a differential graded associative algebra with an
invariant and symmetric co-inner product is an examples of a
V$_\infty$-algebra. In fact, in this case only the lowest maps
$v_2= \mu_2$ and $v_{0,0}=U$ are non zero. That is to say, $v_i$'s
for $i\neq 2$, $v_{(i_1,i_2)}$'s for $(i_1, i_2)\neq (0,0)$, and
$v_{i_1,\cdots,i_k}$'s for $k\geq 3$ are all zero. In this
section, we want to show that in this special case the action of
$\mathcal{DG}_\infty^\bullet$ reduces to the action of the cyclic
Sullivan chord diagrams without marked points, as described in
section \ref{strict-Cyc-section}.

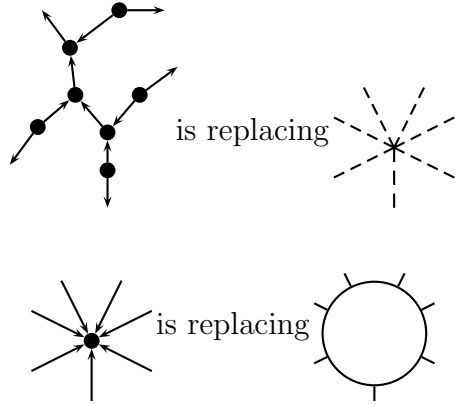
\begin{figure}
\[
     \psset{arrows=<-}
    \pstree[treemode=U, levelsep=0.5]{\Tp}{\pstree[arrows=->]{\Tc*{3pt}}{
    \pstree{\Tc*{3pt}}{ \pstree[treemode=L,arrows=<-]{\Tc*{3pt}}{
    \pstree[treemode=U,levelsep=0.2,arrows=->]{\Tn}{\Tn
    \pstree[levelsep=0.5]{\Tc*{3pt}}{\Tp\pstree[arrows=<-,levelsep=0.2,treemode=R]{\Tn}{
    \pstree[levelsep=0.08]{\Tn}{\pstree[levelsep=0.6,arrows=->]{\Tc*{3pt}}{\Tp}}}}}
    \pstree[levelsep=0.5cm,arrows=->,treemode=D]{\Tc*{3pt}}{\Tp\Tn}}
    \pstree[treemode=R,arrows=<-,levelsep=0]{\Tn}{
    \pstree[levelsep=0.5cm,arrows=->]{\Tc*{3pt}}{\Tp\Tn}} }}}
\begin{pspicture}(0,0)(2,1.6)  
\rput(1,1){is replacing} \end{pspicture}
     \psset{arrows=-}
    \pstree[treemode=U, linestyle=dashed, levelsep=0.8cm, treesep=0.8cm]{\Tp}{
    \pstree[levelsep=0]{\Tp}{
    \pstree[levelsep=0.08cm,treemode=R]{\Tn}{\Tp\Tp}
    \pstree[levelsep=0.8cm,treemode=U]{\Tn}{\Tp\Tp}
    \pstree[levelsep=0.08cm,treemode=L]{\Tn}{\Tp\Tp}    }}
 \]

 \[
     \psset{arrows=->}
    \pstree[treemode=U, levelsep=0.8cm, treesep=0.8cm]{\Tp}{
    \pstree[levelsep=0,arrows=<-]{\Tc*{3pt}}{
    \pstree[levelsep=0.08cm,treemode=R]{\Tn}{\Tp\Tp}
    \pstree[levelsep=0.8cm,treemode=U]{\Tn}{\Tp\Tp}
    \pstree[levelsep=0.08cm,treemode=L]{\Tn}{\Tp\Tp}    }}
\begin{pspicture}(0,0)(2,1.6)  
\rput(1,1){is replacing} \end{pspicture}
     \psset{arrows=-}
    \pstree[treemode=U, levelsep=0.9cm, treesep=0.8cm]{\Tp}{
    \pstree[levelsep=0]{\Tc{20pt}}{
    \pstree[levelsep=0.08cm,treemode=R]{\Tn}{\Tp\Tp}
    \pstree[levelsep=0.8cm,treemode=U]{\Tn}{\Tp\Tp}
    \pstree[levelsep=0.08cm,treemode=L]{\Tn}{\Tp\Tp}    }}
\]
\caption{Chords and strings in the strict case} \label{strict}
\end{figure}

Most of the $v_{i_1,\cdots,i_k}$'s vanish, and the only non-zero
graphs in $\mathcal{DG}_\infty^\bullet$ have vertices of types $1$
and $2$ in the lowest degree. Furthermore, note that due to
associativity and bimodule conditions in Figure \ref{associativity}
and Figure \ref{co-inner product}, one can slide type $1$ vertices
along any type $2$ vertices and also along other type $1$ vertices.
By slight abuse of language, we now call the edges together with
vertices of type $1$ and $2$, the chords of the graph. This makes
sense, since by the above remark, the exact order in which type $1$
and $2$ vertices are applied is irrelevant, and a chord is specified
by only knowing the cyclic combinatorics of its attached type $0$
vertices. We will represent chords by dashed lines and vertices of
type $0$ by open circles, as in Figure \ref{strict}.
Non-vanishing graphs in $\mathcal{DG}_\infty^\bullet$ will
therefore consist of circles with chords attached; see Figure
\ref{chord graph}.

Finally, an orientation for a graph reduces to an orientation as
described in section \ref{strict-Cyc-section}, since, for $d=0$,
all the vertices in the graph are of even degree, and shifts in a
chord may be canonically arranged, so that they only apply to the
chords endpoints.

Our graphs now become just as described in section
\ref{strict-Cyc-section}. By the Theorem \ref{V_k-action} from
section \ref{action on CC section}, graphs without marked points
act on the cyclically invariant part of the Hochschild complex.
This implies, for example, a Lie algebra structure on the cyclic
Hochschild cohomology, which was also noticed by L. Menichi in
\cite{M}.

\subsection{Obstructions and operations}\label{obstract-and-operat}

In this section, we want to describe how the higher type vertices
are used to show that the cyclic Hochschild cohomology of a
$V_k$-algebra, $k\geq 3$, is a Lie algebra. Notice that every
$V_2$-algebra $A$ induces an operation $[\cdot ,\cdot
]:CC^\bullet(A) \otimes CC^\bullet (A)\to CC^\bullet(A)$. We
choose the orientation $v\wedge e_1\wedge e_2$ on the vectors
space generated by the (type $2$) vertex and the two edges, so
that reversing the inputs shows that we have a skew-symmetric
operation. This bracket is closed, and therefore descends to the
cyclic Hochschild cohomology $HC^\bullet(A)$. It turns out that
this bracket, in general, does not satisfy Jacobi identity.
However, in the case of a $V_3$-algebra, the vertices of type $3$
provide homotopies for the Jacobi identity.
\begin{prop}
Let $A$ be a $V_3$-algebra and $[\cdot ,\cdot
]:CC^\bullet(A)\otimes CC^\bullet (A)\to CC^\bullet(A)$ denote the
operation labelled by the graph in Figure \ref{bracket}. Then,
$[\cdot ,\cdot ]$ is skew symmetric and satisfies the Jacobi
identity up to homotopy.
\begin{figure}
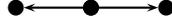

\[
    \psset{arrows=<-}
    \pstree[treemode=L, levelsep=1cm, treesep=0.4cm]{\Tc*{3pt}}{
    \pstree[levelsep=0, arrows=->]{\Tc*{3pt}}{
    \pstree[levelsep=1cm]{\Tc*{3pt}}{\Tc*{3pt}}}}
\]\caption{The bracket $[\cdot,\cdot ]$} \label{bracket}
\end{figure}
\end{prop}
\begin{proof}
In order to see that the Jacobi identity holds, we need to compare
$[[\cdot ,\cdot ],\cdot ]$ with its two cyclic rotations. The
graphs for the term $[[\cdot ,\cdot ],\cdot ]$ of the Jacobi
identity are shown in Figure \ref{jacobi}.
\begin{figure}
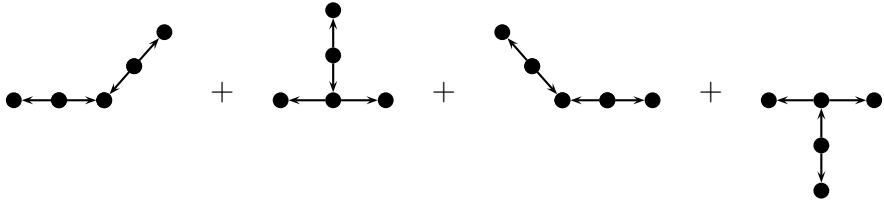

\[
    \psset{arrows=<-}
    \pstree[treemode=R,levelsep=0.6cm, treesep=0.4cm]{\Tc*{3pt}}{
    \pstree[levelsep=0,arrows=->]{\Tc*{3pt}}{\pstree[levelsep=0.6cm]{\Tc*{3pt}}{
    \pstree[levelsep=0,arrows=<-]{\Tc*{3pt}}{\pstree[levelsep=0.4cm]{\Tc*{3pt}}
    {\pstree[levelsep=0,arrows=->]{\Tc*{3pt}}{\pstree[levelsep=0.4cm]{\Tc*{3pt}}
    {\Tc*{3pt}\Tn\Tn}}\Tn\Tn}} }}}
\quad + \quad
    \psset{arrows=->}
    \pstree[treemode=U, levelsep=0, treesep=0.4cm]{\Tc*{3pt}}{
    \pstree[treemode=L, levelsep=0.3cm]{\Tn}{\Tc*{3pt}}
    \pstree[levelsep=0.6cm, arrows=<-]{\Tn}{
    \pstree[arrows=->]{\Tc*{3pt}}{\Tc*{3pt}}}
    \pstree[treemode=R, levelsep=0.3cm]{\Tn}{\Tc*{3pt}}}
\quad + \quad
    \psset{arrows=<-}
    \pstree[treemode=L,levelsep=0.6cm, treesep=0.4cm]{\Tc*{3pt}}{
    \pstree[levelsep=0,arrows=->]{\Tc*{3pt}}{\pstree[levelsep=0.6cm]{\Tc*{3pt}}{
    \pstree[levelsep=0,arrows=<-]{\Tc*{3pt}}{\pstree[levelsep=0.4cm]{\Tc*{3pt}}
    {\pstree[levelsep=0,arrows=->]{\Tc*{3pt}}{\pstree[levelsep=0.4cm]{\Tc*{3pt}}
    {\Tc*{3pt}\Tn\Tn}}\Tn\Tn}} }}}
\quad + \quad
    \psset{arrows=->}
    \pstree[treemode=D, levelsep=0, treesep=0.4cm]{\Tc*{3pt}}{
    \pstree[treemode=L, levelsep=0.3cm]{\Tn}{\Tc*{3pt}}
    \pstree[levelsep=0.6cm, arrows=<-]{\Tn}{
    \pstree[arrows=->]{\Tc*{3pt}}{\Tc*{3pt}}}
    \pstree[treemode=R, levelsep=0.3cm]{\Tn}{\Tc*{3pt}}}
\]
\caption{The term $[[\cdot ,\cdot ],\cdot ]$ for the Jacobi
identity} \label{jacobi}
\end{figure}
The first and third terms only differ by their enumeration of
their three inputs. As required by the Jacobi identity, we also
obtain graphs where these enumerations are cyclically rotated. A
comparison of their orientations shows the cancellation of these
six graphs in the Jacobi identity.

The remaining second and fourth term of $[[\cdot ,\cdot ],\cdot ]$
together with their cyclically rotated input enumerations are
exactly the terms of the boundary of the type $3$ vertex; see
Figure \ref{type-3-boundary-}.
\begin{figure}
\[
\begin{pspicture}(0,0)(10,7)
 \psline(5,5)(3,2) \psline(5,5)(7,2) \psline(3,2)(7,2)
 \psdots[dotsize=6pt](5,3)  \psdots[dotsize=6pt](5,3.7)
 \psdots[dotsize=6pt](4.5,2.5)  \psdots[dotsize=6pt](5.5,2.5)
 \psline[arrows=->](5,3)(5,3.6)     \rput(5.3,3.8){1}
 \psline[arrows=->](5,3)(4.54,2.54) \rput(5.9,2.4){2}
 \psline[arrows=->](5,3)(5.46,2.54) \rput(4.1,2.4){3}
 \psdots[dotsize=6pt](5,1)     \psdots[dotsize=6pt](5,1.6)
 \psdots[dotsize=6pt](5,0.4)   \psdots[dotsize=6pt](5.6,0.4)
 \psdots[dotsize=6pt](4.4,0.4)
 \psline[arrows=->](5,1)(5,1.5)     \psline[arrows=->](5,1)(5,0.5)
 \psline[arrows=->](5,0.4)(5.5,0.4) \psline[arrows=->](5,0.4)(4.5,0.4)
 \rput(5.3,1.6){1} \rput(5.8,0.8){2} \rput(4.2,0.8){3}
 \psdots[dotsize=6pt](7,5)     \psdots[dotsize=6pt](7,4)
 \psdots[dotsize=6pt](8,4)   \psdots[dotsize=6pt](7.5,4.5)
 \psdots[dotsize=6pt](6.5,3.5)
 \psline[arrows=->](7,4)(6.56,3.56)     \psline[arrows=->](7,4)(7.44,4.44)
 \psline[arrows=->](7.5,4.5)(7.06,4.94) \psline[arrows=->](7.5,4.5)(7.94,4.06)
 \rput(6.6,4.8){1} \rput(8,3.6){2} \rput(6.8,3.4){3}
 \psdots[dotsize=6pt](3,5)     \psdots[dotsize=6pt](3,4)
 \psdots[dotsize=6pt](2,4)   \psdots[dotsize=6pt](2.5,4.5)
 \psdots[dotsize=6pt](3.5,3.5)
 \psline[arrows=->](3,4)(3.44,3.56)     \psline[arrows=->](3,4)(2.56,4.44)
 \psline[arrows=->](2.5,4.5)(2.94,4.94) \psline[arrows=->](2.5,4.5)(2.06,4.06)
 \rput(3.4,4.8){1} \rput(2,3.6){3} \rput(3.2,3.4){2}
 \psdots[dotsize=4pt](2,2)     \psdots[dotsize=4pt](2,1)
 \psdots[dotsize=4pt](3,1)     \psdots[dotsize=4pt](2,1.5)
 \psdots[dotsize=4pt](2.5,1)   \psdots[dotsize=4pt](1.6,0.6)
 \psline[arrows=->](2,1.5)(2,1.1) \psline[arrows=->](2,1.5)(2,1.9)
 \psline[arrows=->](2.5,1)(2.1,1) \psline[arrows=->](2.5,1)(2.9,1)
 \psline[arrows=->](2,1)(1.66,0.66)
 \rput(1.7,2){1} \rput(2.9,0.7){2} \rput(1.3,0.8){3}
 \psdots[dotsize=4pt](8,2)     \psdots[dotsize=4pt](8,1)
 \psdots[dotsize=4pt](7,1)     \psdots[dotsize=4pt](8,1.5)
 \psdots[dotsize=4pt](7.5,1)   \psdots[dotsize=4pt](8.4,0.6)
 \psline[arrows=->](8,1.5)(8,1.1) \psline[arrows=->](8,1.5)(8,1.9)
 \psline[arrows=->](7.5,1)(7.9,1) \psline[arrows=->](7.5,1)(7.1,1)
 \psline[arrows=->](8,1)(8.34,0.66)
 \rput(8.2,2){1} \rput(7.1,0.7){3} \rput(8.7,0.8){2}
 \psdots[dotsize=4pt](5,6.9)     \psdots[dotsize=4pt](5,6.4)
 \psdots[dotsize=4pt](5.4,6)     \psdots[dotsize=4pt](4.6,6)
 \psdots[dotsize=4pt](4.2,5.6)   \psdots[dotsize=4pt](5.8,5.6)
 \psline[arrows=->](5.4,6)(5.06,6.34) \psline[arrows=->](4.6,6)(4.94,6.34)
 \psline[arrows=->](5.4,6)(5.74,5.66) \psline[arrows=->](4.6,6)(4.26,5.66)
 \psline[arrows=->](5,6.4)(5,6.8)
 \rput(5.3,6.8){1} \rput(6,5.9){2} \rput(4,5.9){3}
\end{pspicture}
\] \caption{The homotopy induced by the vertex of type $3$}
\label{type-3-boundary-}\end{figure}
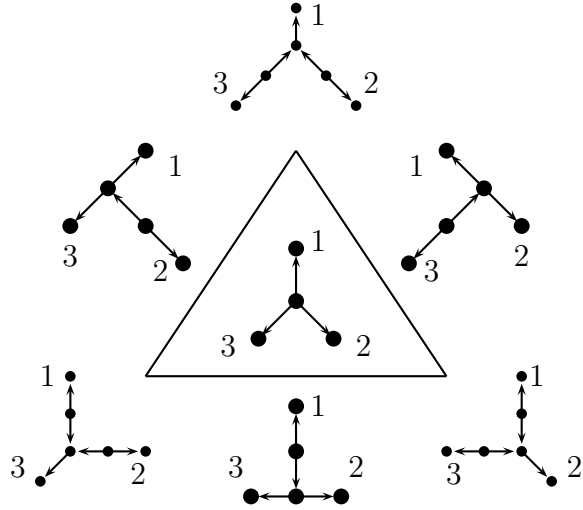
\end{proof}

\section{Some variations}\label{variations section}

We now consider several variations of the discussion of the
previous section. In particular, we will consider loops
(Hochschild complex) and open strings (two-sided cobar complex).
Finally, in section \ref{topological section}, we give a
topological picture which might lend insight to some of our
constructions.

\subsection{Open/closed string operations for V$_\infty$-algebras}\label{open section}

So far, we have only considered the case of the cyclic Hochschild
complex, or the closed strings. It is natural to extend this to
the cyclic Hochschild complex coupled with a two-sided cobar
complex, or an open/closed version of the string topology. The
interaction of open and closed strings is of considerable interest
(see for example, \cite{BCR}, \cite{C}, \cite{H}, \cite{KaSt}, or
\cite{Su}).

\begin{figure}
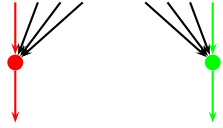

\[
 \psset{arrows=<-, linecolor=red}
    \pstree[treemode=U, levelsep=0.8cm, treesep=0.3cm]{\Tp}
    {\pstree[linecolor=black]{\Tc*{3pt}}{\Tn\Tn\Tn
    \pstree[levelsep=0, linecolor=red]{\Tn} {\Tp}\Tp\Tp\Tp} }
\quad\quad
 \psset{arrows=<-, linecolor=green}
    \pstree[treemode=U, levelsep=0.8cm, treesep=0.3cm]{\Tp}
    {\pstree[linecolor=black]{\Tc*{3pt}}{\Tp\Tp\Tp
    \pstree[levelsep=0, linecolor=green]{\Tn} {\Tp}\Tn\Tn\Tn} }
\] \caption{Homotopy left and right module structures}
\label{homotopy left and right module}
\end{figure}

We now show how to extend our discussion to the open case. Let us
start with the generalizations of the concepts of left and right
modules. Let $A$ be a differential graded associative algebra, and
let $M$ and $N$ be differential graded left and right modules,
respectively. In other words, there are chain maps $M\otimes A\to
M$, and $A\otimes N\to N$, respecting the product in $A$. The
homotopy versions require sequences of maps $\{\lambda_n: M\otimes
A^{\otimes n}\to M\}_{n\geq 0}$, and $\{\rho_n:A^{\otimes
n}\otimes N \to N\}_{n\geq 0}$, where $\lambda_0$ and $\rho_0$ are
the differentials $\partial_M$ and $\partial_N$ of $M$ and $N$,
respectively. Graphically, we will distinguish $A$, $M$ and $N$ by
assigning the colors black, red, and green to $A$, $M$ and $N$,
respectively. Then, the left and right module structures are shown
as in Figure \ref{homotopy left and right module},
and the required relations are as indicated in Figure \ref{homotopy
left and right relations}.
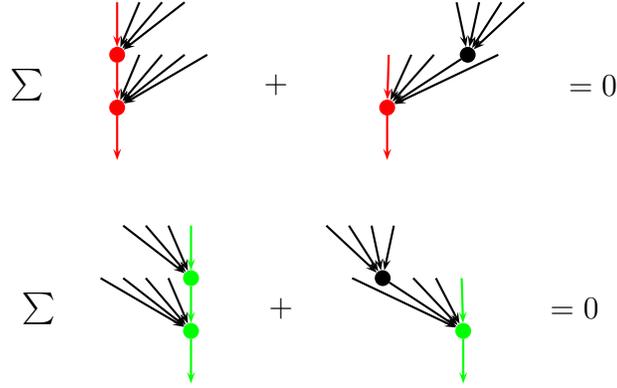
\begin{figure}
\[
\begin{pspicture}(0,0)(.5,2) \rput(.5,1){$\sum$} \end{pspicture}
    \psset{arrows=<-, linecolor=red}
    \pstree[treemode=U, levelsep=0.7cm, treesep=0.3cm]{\Tp}
    {\pstree[linecolor=black]{\Tc*{3pt}}{\Tn\Tn\Tn\Tn
    \pstree[levelsep=0, linecolor=red]{\Tn} {
    \pstree[levelsep=0.7cm,linecolor=black]{\Tc*{3pt}}{\Tn\Tn\Tn
       \pstree[levelsep=0,linecolor=red]{\Tn}{\Tp}\Tp\Tp\Tp}}\Tp\Tp\Tp\Tp} }
\quad\begin{pspicture}(0,0)(.5,2) \rput(.5,1){$+$} \end{pspicture}
 \psset{arrows=<-, linecolor=red}
    \pstree[treemode=U, levelsep=0.7cm, treesep=0.3cm]{\Tp}
    {\pstree[linecolor=black]{\Tc*{3pt}}{\Tn\Tn\Tn\Tn\Tn
    \pstree[levelsep=0, linecolor=red]{\Tn} {\Tp}\Tp\Tp\pstree{\Tc*{3pt}}{\Tn\Tn\Tp\Tp\Tp\Tp}\Tp} }
\quad\begin{pspicture}(0,0)(.5,2) \rput(.5,1){$=0$}\end{pspicture}
\]

\[
\quad\begin{pspicture}(0,0)(.5,2) \rput(.5,1){$\sum$}
\end{pspicture} \quad\quad \psset{arrows=<-, linecolor=green}
    \pstree[treemode=U, levelsep=0.7cm, treesep=0.3cm]{\Tp}
    {\pstree[linecolor=black]{\Tc*{3pt}}{\Tp\Tp\Tp\Tp
    \pstree[levelsep=0, linecolor=green]{\Tn} {
    \pstree[levelsep=0.7cm,linecolor=black]{\Tc*{3pt}}{\Tp\Tp\Tp
       \pstree[levelsep=0,linecolor=green]{\Tn}{\Tp}\Tn\Tn\Tn}}\Tn\Tn\Tn\Tn} }
\begin{pspicture}(0,0)(.6,2) \rput(0,1){$+$} \end{pspicture}
 \psset{arrows=<-, linecolor=green}
    \pstree[treemode=U, levelsep=0.7cm, treesep=0.3cm]{\Tp}
    {\pstree[linecolor=black]{\Tc*{3pt}}{\Tp\pstree{\Tc*{3pt}}{\Tp\Tp\Tp\Tp\Tn\Tn}\Tp\Tp
    \pstree[levelsep=0, linecolor=green]{\Tn} {\Tp}\Tn\Tn\Tn\Tn\Tn} }
    \begin{pspicture}(0,0)(.5,2) \rput(0,1){$=0$} \end{pspicture}
\] \caption{Homotopy left and right module relations} \label{homotopy
left and right relations}
\end{figure}

Next, look at the two-sided cobar complex, $$ C^\bullet(M,A,N)=
\prod_{j=0}^\infty \Big\{f:M\otimes A^{\otimes j}\otimes N\to {\bf
k}\Big\}.$$ The differential $\delta:C^\bullet(M,A,N)\to
C^\bullet(M,A,N)$ is given by,
\begin{eqnarray*}
\delta f(m, a_1,\cdots,a_j,n)&=&\quad \sum_{l}\pm f\big(
\lambda_l(m,a_1,\cdots, a_l),\cdots, a_j,n \big)\\&&+ \sum_{l,i}\pm
f \big(m,a_1, \cdots, \mu_l(a_i, \cdots, a_{i+l-1}), \cdots , a_j,n
\big)\\&&+\sum_{l}\pm f\big (m,a_1,\cdots, \rho_l (a_{j-l+1},\cdots,
a_j,n) \big).
\end{eqnarray*}
Once again, the signs are again given by comparing the linear
order of the symbols in an expression to the orientation $\psi
\wedge f \wedge m \wedge a_1 \wedge \cdots \wedge a_{j} \wedge n$,
where $\psi =\mu_l$, $\lambda_l$, or $\rho_l$, is of degree $1$,
and $m$, $n$, as well as the $a_i$'s, have shifted degrees. The
elements of $C^\bullet(M,A,N)$ will be represented by vertices
with $j$ black incoming edges together with one red and one green
incoming edge on the left and right, respectively; see Figure
\ref{cobar complex}.

\begin{figure}
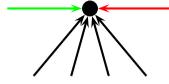

\[
    \psset{arrows=->}
    \pstree[treemode=U, levelsep=0.9cm, treesep=0.5cm]{\Tn}{
    \pstree[levelsep=0]{\Tc*{3pt}}{
    \pstree[levelsep=0.6cm,treemode=L,arrows=<-,linecolor=green]{\Tn}{\Tp}
    \pstree[levelsep=0.9cm,treemode=D,arrows=<-]{\Tn}{\Tp\Tp\Tp\Tp}
    \pstree[levelsep=0.6cm,treemode=R,arrows=<-,linecolor=red]{\Tn}{\Tp}
    }}
 \] \caption{The cobar complex $C^\bullet(M,A,N)$} \label{cobar
complex}
\end{figure}
\begin{figure}
\[
\begin{pspicture}(0,0)(5.2,3.4)
 \psdots[dotsize=6pt,linecolor=green](0.4,3.2) \psdots[dotsize=6pt,linecolor=green](1.4,3.2)
 \psdots[dotsize=6pt](2.4,3.2) \psdots[dotsize=6pt,linecolor=red](3.4,3.2)
 \psdots[dotsize=6pt,linecolor=red](0.8,0.2) \psdots[dotsize=6pt,linecolor=red](1.8,0.2)
 \psdots[dotsize=6pt](2.8,0.2) \psdots[dotsize=6pt,linecolor=green](3.8,0.2)
 \psdots[dotsize=6pt](1.8,1.2) \psdots[dotsize=6pt](2.8,1.2)
 \psdots[dotsize=6pt](1.8,2.2) \psdots[dotsize=6pt](2.8,2.2)
 \psdots[dotsize=6pt](3.8,1.8) \psdots[dotsize=6pt](4.8,1.6)
 \psdots[dotsize=6pt,linecolor=red](4.6,0.7) \psdots[dotsize=6pt,linecolor=green](5.0,2.5)
 \psline[arrows=->,linecolor=green](0.4,3.2)(1.3,3.2)
 \psline[arrows=->,linecolor=green](1.4,3.2)(2.3,3.2)
 \psline[arrows=->,linecolor=red](3.4,3.2)(2.5,3.2)
 \psline[arrows=->,linecolor=red](0.8,0.2)(1.7,0.2)
 \psline[arrows=->,linecolor=red](1.8,0.2)(2.7,0.2)
 \psline[arrows=->,linecolor=green](3.8,0.2)(2.9,0.2)
 \psline[arrows=->,linecolor=red](4.6,0.7)(4.78,1.52)
 \psline[arrows=->,linecolor=green](5.0,2.5)(4.82,1.68)
 \psline[arrows=->](3.8,1.8)(4.7,1.6) \psline[arrows=->](3.8,1.8)(2.9,2.2)
 \psline[arrows=->](2.8,2.2)(2.5,3.17) \psline[arrows=->](2.8,2.2)(2.8,1.3)
 \psline[arrows=->](1.8,1.2)(1.8,0.3) \psline[arrows=->](2.8,1.2)(2.8,0.3)
 \psline[arrows=->](1.8,1.2)(2.7,2.2) \psline[arrows=->](2.8,1.2)(1.8,2.1)
 \psline[arrows=->](1.8,2.2)(1.4,3.1) \psline[arrows=->](1.8,2.2)(2.4,3.1)
\end{pspicture}
\]
\caption{An open/closed graph} \label{typical open closed}
\end{figure}
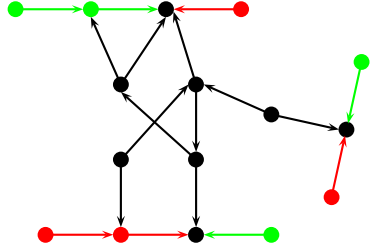

We want to consider all graphs, consisting of vertices of the
above type. We also need to add red and green edges with only one
outgoing edge and no incoming edges, corresponding to the identity
map $M\to M$ or $N\to N$, respectively. In this manner, we obtain
graphs similar to that of Figure \ref{typical open closed}.

Let $\coDG$ denote the space of graphs that include the above
colored vertices, with orientations given by all vertices and all
edges of type $n\neq 0$, where the red and green vertices from
Figure \ref{homotopy left and right module} have degree $-2$.
Furthermore, $\coDG$ has a differential and a composition, which
are defined similarly to those in section \ref{graph complex
section}. With these adjustments, we can extend the theorem from
section \ref{action on CC section} in the following way.
\begin{thrm}\label{open-closed-theorem}
There is a map,
\begin{multline*}
\coDG \to \prod_{r_c,r_o,s_c,s_o}Hom \Big( CC^ \bullet(A)^{ \otimes
r_c}\otimes C^\bullet(M,A,N)^{\otimes r_o},\\ CC^\bullet(A)^{\otimes
s_c}\otimes C^\bullet(M,A,N)^{\otimes s_o}\Big),
\end{multline*}
which respects the grading, differential, and composition.
\end{thrm}

\begin{remk}
The above theorem implies that $H^\bullet(M,A,N)$, the homology of
the complex $C^\bullet(M, A, N)$, is a module over
$HC^\bullet(A)$. In section \ref{topological section}, we will see
how $M$ and $N$ correspond to subcomplexes of a given manifold;
see Figure \ref{topo-open-closed}. The module structure described
in this section corresponds to the intersection of cells inside
the big manifold.

But, this is only half the story. When the subcomplexes are in
fact submanifolds, then they will also have Poincar\'e duality,
making it possible for intersection within them. Algebraically, we
might add this structure by requiring $M$ and $N$ to also have
(colored) homotopy co-inner products. This would give a further
extension of the PROP $\coDG$, with a similar action as in Theorem
\ref{open-closed-theorem}.
\end{remk}

\subsection{Algebraic structure of the Hochschild complex}\label{non-cyclic section}

In this section, we consider the Hochschild complex $CH^\bullet(A,
A^\ast)$. A modification of the string PROP, that corresponds to
introducing marked points, will give the desired PROP action in
the non-cyclic setting.

Let $A$ be a V$_\infty$-algebra. An element $\1\in A$ is called a
strict unit, or simply a unit, if it satisfies the relations,
\begin{eqnarray*}
 v_2(a,\1)=v_2(\1,a)=a && \text{for } a\in A,\\
v_{i_1,\cdots,i_k}(\cdots,\1,\cdots)=0 && \text{for all other }
(i_1,\cdots,i_k)\neq (2).
\end{eqnarray*}
The Hochschild cochain complex of $A$ with values in the dual space
$A^*$ is given by,
$${CH^\bullet}(A,A^*)=\prod_{j=0}^\infty \Big\{f:A^{ \otimes
j}\to A^\ast \Big\}.$$ In the presence of a unit, the normalized
subcomplex is defined by,
\begin{multline*}
\overline{CH^\bullet}(A,A^*)=\prod_{j=0}^\infty \Big\{f:A^{
\otimes j}\to A^\ast \Big|
\quad\quad\quad\quad\quad\quad\quad\quad\quad
\quad\quad\quad\quad\,\, \\%
f(a_1, \cdots, a_j)=0, \text{ if any of the } a_1,\cdots, a_j
\text{ equals }\1 \Big\}.
\end{multline*}
The induced differential on this subspace makes it into a
subcomplex, which is quasi-isomorphic to the big complex via the
inclusion $\overline{CH^\bullet}(A,A^*) \hookrightarrow
CH^\bullet(A,A^*)$; see \cite{L}. In what follows, we identify a
map $A^{\otimes j} \to A^\ast$ with a map $A^{\otimes j+1}\to {\bf
k}$.

We now define a new graph complex, $\mDG$, suitable for treating
the normalized Hochschild complex. $\mDG$ is a variation of
$\mathcal{DG}_\infty^\bullet$ that deals with the issues of
starting points for inputs and outputs, as well as the unit. Here
are what we need,

\begin{enumerate}
\item We need marked point at every vertex of type $0$. By
definition, a marked point of an input is a choice of edge
attached to that input. For the trivial type $0$ vertex without
edges, no choice is necessary. We will denote this chosen edge by
a double-headed arrow; see Figure \ref{marked-DG}. \item A marked
point for an output consists of exactly one additional external
leg for each output. We denote this leg by a vertex with one
outgoing edge and no incoming edges; see Figure
\ref{output-marked-points}.
\begin{figure}
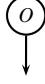

\[
\psset{arrows=->}
\pstree[treemode=D,levelsep=0.8cm]{\Tcircle{o}}{\Tp}
\] \caption{The output marked point} \label{output-marked-points}
\end{figure}
\item  The unit $\1\in A$ gives rise to a new vertex with exactly
one outgoing edge and no incoming edges; see Figure
\ref{the-unit}.
\begin{figure}
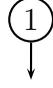

\[
\psset{arrows=->}
\pstree[treemode=D,levelsep=0.8cm]{\Tcircle{1}}{\Tp}
\] \caption{The unit vertex} \label{the-unit}
\end{figure}
By the defining relations of the unit, we may never attach the unit
to any $v_{i_1, \cdots, i_k}$, except for $v_{i}$, when $i=2$, for
which we may ignore the unit, since $v_2(a,\1)=v_2(\1,a)=a$. Thus,
the unit may only be applied to vertices of type $0$. Furthermore,
since we work with the normalized Hochschild complex, the unit can
only be applied in the last slot $f(\cdots ,\1)$ of any $f\in
\overline{ HC^\bullet}(A,A^*)$. This means that if the unit vertex
occurs, it necessarily has to be the marked point of a type $0$
vertex.
\end{enumerate}
An example of a graph in $\mDG$ is shown in Figure
\ref{marked-DG}.
\begin{figure}
\[
\begin{pspicture}(2.5,0.5)(8.4,4)
%
 \psdots[dotsize=6pt](3,3)     \psdots[dotsize=6pt](7.2,2)
 \psdots[dotsize=6pt](4.5,4)   \psdots[dotsize=6pt](4,3)
 \psdots[dotsize=6pt](5.5,0.5) \psdots[dotsize=6pt](5.2,1.5)
 \psdots[dotsize=6pt](5,3)     \psdots[dotsize=6pt](6.5,3)
 \psdots[dotsize=6pt](7.6,3.8) \psdots[dotsize=6pt](4.5,2)
 \rput(3,2){$o$}     \rput(6.5,2){$o$}     \rput(7.9,2.8){$1$}
 \pscircle(3,2){0.25}\pscircle(6.5,2){0.25}\pscircle(7.9,2.8){0.25}
 \psline[arrows=->](4,3)(4.9,3)
 \psline[arrows=->](6.5,3)(5.1,3)
 \psline[arrows=->](5.2,1.5)(5.44,0.56)
 \psline[arrows=->](5.2,1.5)(5,2.9)
 \psline[arrows=->](4.5,2)(4.5,3.9)
 \psline[arrows=->>](4,3)(3.1,3)
 \psline[arrows=->](3,2.25)(3,2.9)
 \psline[arrows=->](6.5,2.25)(6.5,2.9)
 \pscurve[arrows=->](7.2,2)(7,2.8)(6.6,3)
 \pscurve[arrows=->](7.2,2)(6.8,1)(5.6,0.5)
 \pscurve[arrows=->](6.5,3)(5.6,3.7)(4.6,4)
 \pscurve[arrows=->](4.5,4)(3.4,3.6)(3.04,3.06)
 \pscurve[arrows=->](4.5,2)(4.7,0.7)(5.4,0.5)
 \pscurve[arrows=->](6.5,3)(6.9,3.6)(7.5,3.8)
 \pscurve[arrows=->](5,3)(4.7,3.6)(4.53,3.93)
 \pscurve[arrows=->>](4.5,4)(5.3,3.6)(5.9,1.6)(5.56,0.56)
 \pscurve[arrows=->>](7.9,3.05)(7.8,3.4)(7.66,3.74)
\end{pspicture}
\] \caption{A graph in $\mDG$} \label{marked-DG}
\end{figure}
The differential on $\mDG$ is again given by expansion of an edge,
where we now have to preserve the marked points of the inputs and
outputs. Similarly, when composing the graphs $\gamma_1$ and
$\gamma_2$, we also need to match the marked edges. This means
that when attaching the edges from some input of $\gamma_2$ to the
corresponding output of $\gamma_1$, the marked edge of the input
has to coincide with that of the output. With this, we have the
following theorem.
\begin{thrm}
There is a map,
\begin{equation*}
\mDG \to \prod_{r,s} Hom\Big(\overline{HC^\bullet} (A,A^*)^{\otimes
r},\overline{HC^\bullet}(A,A^*)^{\otimes s}\Big),
\end{equation*}
which respects the grading, differential, and composition.
\end{thrm}

The following comments are in order.
\begin{remk}
Let $\1\in A$ denote the unit. We can extend the definition of a
V$_\infty$-algebra $\{v_{i_1,\cdots,i_k}\}_k$ satisfying all the
relations from Definition \ref{V_k}, by allowing for $v_0$ to
exist; see Remark \ref{weak-V-infty}. Here, the element $v_0\in A$
will be called a weak unit which makes $A$ into a weak
A$_\infty$-algebra.
\end{remk}

\begin{remk}
In the case of strictly associative algebras, with the help of the
enumeration of the inputs and the linear ordering of the marked
points on the input circles, the orientation on vertices and edges
may be chosen canonically. Thus, in this case, our constructions
and results reduce to those of section \ref{strict-Hoch-section}.
\end{remk}

\begin{remk}\label{kont-soib-M}
The graph complex $\mathcal{DG}_\infty^\bullet$ is reminiscent of
that considered in \cite{KoSo} by Kontsevich and Soibelman. They
define an operad $M$ that naturally acts on the Hochschild cochain
complex $ CH^\bullet(A,A)= \prod_{j=0}^\infty \big\{f:A^{\otimes
j}\to A\big\}$ of an A$_\infty$-algebra $A$. The operad $M$ is
made out of trees that have only type $1$ vertices which are used
for both the components of the A$_\infty$-algebra structure as
well as for labelling the elements of $CH^\bullet(A,A)$.

Their discussion in \cite{KoSo} goes further by identifying $M$
with the Strebel differentials on the Riemann sphere. Thus, a
quasi-isomorphism between the operad $M$ and the chains on the
little disc operad is established, which in turn solves the
Deligne conjecture. We suspect a strong relationship between
$\mathcal{DG}_\infty^\bullet$ and the chains on a moduli space of
Riemann surfaces.
\end{remk}

\begin{remk}
This work is about the associative operad. It is reasonable to
expect that similar results would hold for any cyclic operad
$\mathcal O$. More precisely, we expect that the space of directed
graphs with an $\mathcal O$-cyclic order on their vertices would
act on the cyclic cochain complex associated to the operad
$\mathcal O$. See \cite{LT} for a study of homotopy co-inner
products on cyclic operads.
\end{remk}

\subsection{Topological applications and motivations}\label{topological section}
Let $X$ denote a simply connected and triangulated homology
manifold of dimension $d$. In a previous work \cite{Z}, we showed
that the (simplicial) cochains $C^\bullet X$ form a V$_2$-algebra,
i.e., an A$_\infty$-algebra with an invariant and symmetric
homotopy co-inner product. This was achieved by obtaining an
invariant and symmetric homotopy inner product that was moreover
nondegenerate in an appropriate sense. It was then argued, using
minimal models, how this nondegeneracy gives rise to an invariant
and symmetric homotopy co-inner product. Note that in this paper
co-inner products do not satisfy any nondegeneracy assumptions. An
immediate application of Theorem \ref{V_k-action} is the following
corollary.
\begin{coro}
The PROP $\mathcal{DG}_2$ acts on the cyclic Hochschild complex
$CC^\bullet(C^\bullet X)$.
\end{coro}

\begin{coro}
The PROP ${\it m}\mathcal{DG}_2$ acts on the Hochschild complex
$CH^\bullet(C^\bullet X)$.
\end{coro}

It is known that for a simply connected $X$, the Hochschild
complex $CH^\bullet( C^\bullet X, C_\bullet X)$ and the cyclic
Hochschild complex $CC^\bullet(C^\bullet X)$ are quasi-isomorphic
to the chains on the free loop space, $C_\bullet LX$, and the
equivariant chains on the free loop space, $C_\bullet^{S^1}LX$,
respectively; see Figure \ref{topo-marked} and Figure
\ref{topo-closed}. Therefore, we have,
\begin{figure}
\[
\begin{pspicture}(0.2,0)(3,2.4)
 \pscurve(0.4,0.2)(1.5,0.3)(2.6,0)  \pscurve(2.6,0)(2.6,1.4)(3,2)
 \pscurve(0.8,2.2)(1.9,2.3)(3,2)  \pscurve(0.4,0.2)(0.5,1.6)(0.8,2.2)
 \psellipse(1.55,1.3)(0.45,0.35) \psline(1.5,1.05)(1.5,0.87)
\end{pspicture}
\] \caption{A loop in a manifold} \label{topo-marked}
\end{figure}
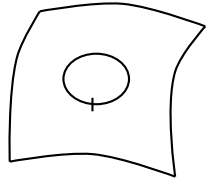
\begin{figure}
\[
\begin{pspicture}(0.2,0)(3,2.4)
 \pscurve(0.4,0.2)(1.5,0.3)(2.6,0)  \pscurve(2.6,0)(2.6,1.4)(3,2)
 \pscurve(0.8,2.2)(1.9,2.3)(3,2)  \pscurve(0.4,0.2)(0.5,1.6)(0.8,2.2)
 \psellipse(1.55,1.3)(0.45,0.35)
\end{pspicture}
\] \caption{A closed string in a manifold} \label{topo-closed}
\end{figure}
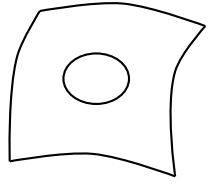

\begin{coro}
The PROP $\mathcal{DG}_2$ acts on $C_\bullet^{S^1}LX$.
\end{coro}
\begin{coro}
The PROP ${\it m}\mathcal{DG}_2$ acts on $C_\bullet LX$.
\end{coro}
\begin{coro}
$H_\bullet LX$ is, after an appropriate shift in degrees, a BV
algebra.
\end{coro}

\begin{figure}
\[
\begin{pspicture}(0.2,0)(3,2.4)
 \pscurve(0.4,0.2)(1.5,0.3)(2.6,0)  \pscurve(2.6,0)(2.6,1.4)(3,2)
 \pscurve(0.8,2.2)(1.9,2.3)(3,2)  \pscurve(0.4,0.2)(0.5,1.6)(0.8,2.2)
 \psellipse(1.55,1.5)(0.3,0.2) \pscurve(0.8,0.8)(1.6,0.83)(2.1,0.7)
 \pscurve[linecolor=red](0.8,0.4)(0.9,1.4)(1.2,2)
 \pscurve[linecolor=green](2.1,0.3)(2.2,1.3)(2.5,1.9)
\end{pspicture}
\] \caption{Open/closed strings in a manifold} \label{topo-open-closed}
\end{figure}
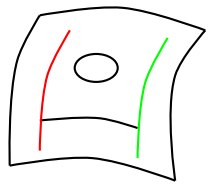

Also, for simplicial subcomplexes $K, L \subset X$, let $P(K,X,L)$
be the space of paths in $X$ starting in $K$ and ending in $L$.
The cochains $C^\bullet K$ and $C^\bullet L$ are left and right
modules over $C^\bullet X$, and the two-sided cobar complex
$C^\bullet( C^\bullet K, C^\bullet X, C^\bullet L)$ is naturally
quasi-isomorphic to the chains of the path space $C_\bullet (P(K,
X ,L))$; see Figure \ref{topo-open-closed}.

We believe there is a deeper connection between string topology
and the structures discussed in this paper. More specifically, we
expect that for a triangulated Poincar\'e duality space $X$, the
simplicial chains $C^\bullet X$ to be a V$_\infty$-algebra, which
in turn would imply that $\mathcal{DG}_\infty^\bullet$, $\mDG$,
and $\coDG$ act on $CC^\bullet(C^\bullet X)$,
$CH^\bullet(C^\bullet X, C_\bullet X)$, and $CC^\bullet(C^\bullet
X) \otimes C^\bullet(C^\bullet K, C^\bullet X, C^\bullet L)$,
respectively.

In light of all this, one also expects an action of
$\mathcal{DG}_\infty^\bullet$ on $C^{S^1}_\bullet LX$, an action
of $\mDG$ on $C_\bullet LX$, and an action of $\coDG$ on $C^{S^1
}_\bullet LX \otimes C_\bullet(P(K, X, L))$. It is highly
desirable to understand the relationship between the PROP
$\mathcal{DG}_\infty^\bullet$ and the PROP of the chains on
Deligne-Mumford compactification of the moduli space of Riemann
surfaces.


\begin{thebibliography}{00000}
\bibitem[BCR]{BCR} N.A. Baas, R.L. Cohen, A. Ramirez, \textit{The
 topology of the category of open and closed strings}, AT/0411080, (2004)
\bibitem[CS1]{CS1} M. Chas, D. Sullivan, \textit{String Topology}, GT/9911159, (1999)
\bibitem[CS2]{CS2} M. Chas, D. Sullivan, \textit{Closed string operators
 in topology leading to Lie bialgebras and higher string algebra},
 The legacy of Niels Henrik Abel,  771--784, Springer, Berlin, (2004)
\bibitem[C]{C} K. Costello, \textit{Topological conformal field
 theories and Calabi-Yau categories}, QA/0412149, (2004)
\bibitem[G]{G} M. Gerstenhaber, \textit{The cohomology structure of an
 associative ring}, Annals of Mathematics, Vol. 78, No. 2, (1963)
\bibitem[H]{H} E. Harrelson, \textit{On the homology of
 open-closed string theory}, AT/0412249, (2004)
\bibitem[Kaj]{K} H. Kajiura, \textit{Noncommutative homotopy algebras
 associated with open strings}, QA/0306332, (2003)
\bibitem[KaSt]{KaSt} H. Kajiura, J. Stasheff, \textit{Homotopy algebras
 inspired by classical open-closed string field theory}, QA/0410291, (2004)
\bibitem[Kau]{Ka1} R.M. Kaufmann, \textit{On spineless cacti, Deligne's
 conjecture and Connes--Kreimer's Hopf algebra}, QA/0308005, (2003)
\bibitem[KoSo]{KoSo} M. Kontsevich, Y. Soibelman, \textit{Deformations of
 algebras over operads and Deligne's conjecture},  Conférence Moshé Flato 1999,
 Vol. I (Dijon),  255--307, Math. Phys. Stud., 21, Kluwer Acad. Publ., Dordrecht, (2000)
\bibitem[L]{L} J.-L. Loday, \textit{Cyclic Homology}, Grundlehren der
 mathematischen Wissenschaften 301, Springer-Verlag, (1992)
\bibitem[LT]{LT} R. Longoni, T. Tradler, \textit{Homotopy inner products for
 cyclic operads}, AT/0312231, (2003)
\bibitem[MS1]{MS1} J.E. McClure, J.H. Smith, \textit{A solution of Deligne's
 Hochschild cohomology conjecture}, Contemp. Math.,Vol 293, 153--193, (2002)
\bibitem[MS2]{MS2} J.E. McClure, J.H. Smith, \textit{Operads
 and cosimplicial objects: an introduction},  Axiomatic, enriched and motivic homotopy theory,
  133--171, NATO Sci. Ser. II Math. Phys. Chem., 131, Kluwer Acad. Publ., Dordrecht, 2004
\bibitem[M]{M} L. Menichi, \textit{Batalin-Vilkovisky algebras and cyclic
 cohomology of Hopf algebras},  $K$-Theory  32  (2004),  no. 3, 231--251
\bibitem[St]{St} J. Stasheff, \textit{On the homotopy associativity of
 H-spaces I and II}, Trans. AMS  108, p. 275-312, (1963)
\bibitem[Su]{Su} D. Sullivan, \textit{Open and closed string field theory
 interpreted in classical algebraic topology},  Topology, geometry and quantum field theory,
 344--357, London Math. Soc. Lecture Note Ser., 308, Cambridge Univ. Press, Cambridge, (2004)
\bibitem[Ta1]{Ta1} D. Tamarkin, \textit{Another proof of M. Kontsevich formality
 theorem}, QA/9803025, (1998)
\bibitem[Ta2]{Ta2} D. Tamarkin, \textit{Formality of chain operad of little discs},
 Lett. Math. Phys.  66 (2003),  no. 1-2, 65--72
\bibitem[Tr1]{T} T. Tradler, \textit{Infinity-inner-products on
 A-infinity algebras}, AT/0108027, (2001)
\bibitem[Tr2]{T2} T. Tradler, \textit{The BV algebra on Hochschild cohomology
 induced by infinity inner products}, QA/0210150, (2002)
\bibitem[TZ1]{Z} T. Tradler, M. Zeinalian, \textit{Poincare duality at the chain
 level, and a BV structure on the homology of the free loop space of
 a simply connected Poincare duality space}, math.AT/0309455,
 (2003)
\bibitem[TZ2]{TZ} T. Tradler, M. Zeinalian, \textit{On the cyclic Deligne
 conjecture}, J. Pure and Appl. Alg., 204 (2006) 280-299
\bibitem[V]{V} A.A. Voronov, \textit{Homotopy Gerstenhaber algebras},
 Conférence Moshé Flato 1999, Vol. II (Dijon),  307--331, Math. Phys. Stud., 22,
 Kluwer Acad. Publ., Dordrecht, (2000)
\end{thebibliography}
\end{document}